\documentclass[12pt]{amsart}

\newtheorem{theorem}{Theorem}[section]
\newtheorem{notations}[theorem]{Notation}
\newtheorem{lemma}[theorem]{Lemma}
\newtheorem{corollary}[theorem]{Corollary}
\theoremstyle{definition}

\theoremstyle{remark}

\numberwithin{equation}{section}
\usepackage[top=1in, bottom=1in, left=1in, right=1in]{geometry}
\usepackage{times}
\usepackage{enumerate}
\newcommand{\colvec}[2][.75]{%
  \scalebox{#1}{%
    \renewcommand{\arraystretch}{.8}%
    $\begin{bmatrix}#2\end{bmatrix}$%
  }
}

\newcommand{\collvec}[2][.8]{%
  \scalebox{#1}{%
    \renewcommand{\arraystretch}{.8}%
    $\begin{bmatrix}#2\end{bmatrix}$%
  }
}
\usepackage{graphicx,adjustbox}
\usepackage{hyperref}
\usepackage{amsmath,amssymb}
\usepackage{amscd}
\usepackage{graphicx}
\usepackage[all]{xy}
\usepackage{xcolor}

\title[Projective closures of affine curves]
{Projective closures of affine monomial curves}
\author{
Joydip Saha
\and
Indranath Sengupta
\and
Pranjal Srivastava
}
\date{}
\address{\small \rm  Stat Math Unit, Indian Statistical Institute, Kolkata, West-Bengal,700108, INDIA.} 
\email{saha.joydip56@gmail.com}
\thanks{The first author thanks NBHM, Government of India for post-doc fellowship at ISI kolkata.}

\address{\small \rm  Discipline of Mathematics, IIT Gandhinagar, Palaj, Gandhinagar, 
Gujarat 382355, INDIA.}
\email{indranathsg@iitgn.ac.in}
\thanks{The second author is the corresponding author.}

\address{\small \rm  Discipline of Mathematics, IIT Gandhinagar, Palaj, Gandhinagar, 
Gujarat 382355, INDIA.}
\email{pranjal.srivastava@iitgn.ac.in}
\date{}



\begin{document}

\begin{abstract}
We study the projective closures 
of three important families of affine monomial curves in 
dimension $4$, namely the Backelin curve, the Bresinsky curve 
and the Arslan curve, 
in order to explore possible connections between 
syzygies and the arithmetic Cohen-Macaulay property.
\end{abstract}

\maketitle

\section{Introduction}

Let $\mathbb{N}$ denote the set of nonnegative integers and $k$ denote a field. Let $r\geq 3$ and $\mathbf{\underline n} = (n_{1}, \ldots, n_{r})$ be a sequence of 
$r$ distinct positive integers with $\gcd(\mathbf{\underline n})=1$. Let us assume 
that the numbers $n_{1}, \ldots, n_{r}$ generate the numerical semigroup 
$\Gamma(n_1,\ldots, n_r) = \lbrace\sum_{j=1}^{r}z_{j}n_{j}\mid z_{j}\in\mathbb{N} \rbrace$ 
minimally, that is, if $n_i=\sum_{j=1}^{r}z_{j}n_{j}$ for some non-negative 
integers $z_{j}$, then $z_{j}=0$ for all $j\neq i$ and $z_{i}=1$. Let 
$\eta:k[x_1,\,\ldots,\, x_r]\rightarrow k[t]$ be the mapping defined by 
$\eta(x_i)=t^{n_i},\,1\leq i\leq r$. Let $\frak{p}(n_1,\ldots, n_r) = \ker (\eta)$. 
Let us assume that $n_{r}>n_{i}$ for all $i<r$, for the sequence $\mathbf{\underline n} = (n_{1}, \ldots, n_{r})$. 
We fix $n_{0}=0$. We define the semigroup $\overline{\Gamma(n_{1},\ldots n_{r})}\subset \mathbb{N}^{2}$, ( $\mathbb{N}$ 
is the set of all non negative integers) generated by $\{(n_{i},n_{r}-n_{i})\mid 0\leq i\leq r\}$. Let us denote 
$\overline{\mathfrak{p}(n_{1},\ldots n_{r})}$ be the kernel of $k$-algebra map 
$\eta^{H}:k[x_{0},\ldots,x_{r}]\longrightarrow k[s,t]$, $\eta^{H}(x_{i})=t^{n_{i}}s^{n_{r}-n_{i}}$. 
Then homogenization of the ideal $\frak{p}(n_1,\ldots, n_r) $ with respect to the variable $x_{0}$ is 
$\overline{\mathfrak{p}(n_{1},\ldots n_{r})}$. Thus the projective curve 
$\{[(a^{n_{r}}:a^{n_{r}-n_{1}}b^{n_{1}}:\cdots:b^{n_{r}})]\in\mathbb{P}^{r}_{k}\mid a,b\in k\}$ 
is the projective closure of the affine curve 
$C(n_{1},\ldots,n_{r}):=\{(b^{n_{1}},\ldots b^{n_{r}})\in \mathbb{A}^{r}_{k}\mid b\in k\}$, 
and we denote it by $\overline{C(n_{1},\ldots,n_{r})}$. We say that the projective curve 
$\overline{C(n_{1},\ldots,n_{r})}$ is arithmetically Cohen-Macaulay if the vanishing ideal 
$\overline{\mathfrak{p}(n_{1},\ldots n_{r})}$ is a Cohen-Macaulay ideal.
\medskip

Let $\beta_{i}(\frak{p}(n_1,\ldots, n_r))$ denotes the $i$-th Betti number of the ideal 
$\frak{p}(n_1,\ldots, n_r)$. Therefore, $\beta_{1}(\frak{p}(n_1,\ldots, n_r))$ denotes the 
minimal number of generators of $\frak{p}(n_1,\ldots, n_r)$. 
For a given $r\geq 3$, let $\beta_{i}(r) = {\rm lub}(\beta_{i}(\frak{p}(n_1,\ldots, n_r))$, 
where the lub is taken over all the sequences of positive integers $n_1,\ldots, n_r$. 
Herzog \cite{herzog} proved that, $\beta_{1}(3)$ is $3$ and it follows easily that 
$\beta_{2}(3)$ is a finite integer as well. Bresinsky \cite{bre} (and \cite{brehoa}) 
defined a class of monomial curves in $\mathbb{A}^{4}$ and proved that $\beta_{1}(4)=\infty$. 
He used this observation to prove that $\beta_{1}(r)=\infty$, for every $r\geq 4$. We 
call this family the \textit{Bresinsky curve}. Subsequently, it has been proved in \cite{mssbetti} 
that all the higher Betti numbers of Bresinsky curve are also not bounded above by a 
fixed integer. F. Arslan introduced the curve $\Gamma_{A_{h}}=\langle h(h+1),h(h+1)+1,(h+1)^{2},(h+1)^{2}+1 \rangle$, 
for $h\geq 2$, in $\mathbb{A}^{4}$, and proved the Cohen-Macaulayness of tangent cone of the curve 
at the origin (see \cite{a}). Subsequent to this, Fr\"{o}berg et al. presented the numerical 
semigroup $\langle s,s+3,s+3n+1,s+3n+2 \rangle$, for $ n \geq 2, r\geq 3n+2, s=r(3n+2)+3$, in 
$\mathbb{A}^{4}$, and proved the unboundedness of the $\rm{type}$ of this numerical semigroup 
(see \cite{fgh}). This class of semigroups was first proposed by Backelin. 
All three families have one strong resemblance that all are monomial curves in 
the affine $4$ space, and there is no upper bound on the minimal generating set 
of the defining ideals of these curves. Recently, Herzog and Stamate \cite{hs} have 
proved that the projective closure of the Bresinsky curve is not arithmetically 
Cohen-Macaulay, but the projective closure of the Arslan curve is arithmetically Cohen-Macaulay, 
using a Gr\"{o}bner basis criterion. In \cite{AA}, Stamate has proved that 
the Betti sequence of the affine Bresinsky curve and affine Arslan curve are 
$(1,4h,8h-4,4h-3)$ and $(1,2h+2,4h,2h-1)$, respectively.
\medskip

This paper is devoted to the study of the projective closures of three families: 
The Backelin curve, the Bresinsky curve and the Arslan curve. We have computed the 
syzygies of the projective closures of all three families and have proved that the 
Betti sequences of the affine curves and their projective closures remain the same 
for the Backelin curve and the Arslan curve, and they differ for the Bresinsky curve. 
In fact, the last Betti number of the projective closure of the Bresinsky curve is $1$. Two main questions 
that form the background of this work are the following: (1) To find a suitable sufficient condition    
on an affine monomial curve so that its projective closure is arithmetically Cohen-Macaulay; 
(2) to find a suitable sufficient condition on an affine monomial curve so that the Betti sequence of the 
affine monomial curve is the same as the Betti sequence of its projective closure. 
This paper is attempt to understand these questions through three very important 
and interesting classes of curves mentioned above. Computations with 
\cite{sing} have helped us understand the structure of the 
syzygies for most of the examples. 
\medskip

This paper has been arranged in the following order: First, we study the affine 
Backelin curve, its Gr\"{o}bner basis and use that to prove that its projective 
closure is arithmetically Cohen-Macaulay, with the help of the Gr\"{o}bner basis 
criterion proved in (Theorem 2.2, \cite{hs}). We then compute the syzygies of the 
projective closure and the Hilbert series of the Backelin curve. Next, we compute 
the syzygies of the projective closure of the Bresinsky curve and finally the 
syzygies of the projective closure of the Arslan curve.

\section{The Backelin Curve and its Gr\"{o}bner Basis}
Let us first begin with a description of Backelin's example of monomial 
curves in the affine space $\mathbb{A}^{4}$. Backelin defined the numerical 
semigroups $\langle s, s+3, s+3n+1, s+3n+2\rangle$, 
for $n\geq 2$, $r\geq 3n+2$ and $s=r(3n+2)+3$. In \cite{fgh}, 
it has been shown that the type of such semigroups are not bounded 
by the embedding dimension. 
\medskip

\begin{notations}{\rm 
We fix some notations. For $n\geq 2$, $r\geq 3n+2$ and $s=r(3n+2)+3$,
\medskip

\noindent $\Gamma_{nr}:=\Gamma(s, s+3, s+3n+1, s+3n+2)$;\\[2mm]
$\overline{\Gamma}_{nr}:=\overline{\Gamma(s, s+3, s+3n+1, s+3n+2)}$;\\[2mm]
$\mathfrak{P}_{nr}:=\frak{p}(s, s+3, s+3n+1, s+3n+2)$;\\[2mm]
$\overline{\mathfrak{P}}_{nr}:=\overline{\frak{p}(s, s+3, s+3n+1, s+3n+2)}$, therefore $\overline{\mathfrak{P}}_{nr}=(\mathfrak{P}_{nr})^{H}$ with respect 
to the variable $x_{0}$;\\[2mm]
$B_{nr}$ denotes the affine Backelin curve defined by the numerical semigroup 
$\Gamma_{nr}$;\\[2mm]
$\overline{B}_{nr}$ denotes the projective closure of the Backelin curves;\\[2mm]
$G(I)$ denotes the unique minimal generating set of a monomial ideal $I$ in 
a polynomial ring over a field. 
}
\end{notations}
 
\medskip
  
\begin{theorem}\textbf{(Gastinger)}\quad \label{gastinger} 
Let $A = k[x_{1},\ldots,x_{r}]$ be the polynomial ring, $I\subset A$ 
the defining ideal of a monomial curve defined by natural numbers 
$a_{1},\ldots,a_{r}$, whose greatest common divisor is $1$.  
Let $J$ be an ideal contained in $I$. Then $J = I$ if and 
only if $\mathrm{dim}_{k} A/\langle J + (x_{i}) \rangle =a_{i}$, 
for some $i$; equivalently 
$\mathrm{dim}_{k} A/\langle J + (x_{i}) \rangle =a_{i}$ for any $i$.
\end{theorem}  

\proof See in \cite{g}.\qed
\medskip
 
Let $I\subset k[x_{1},\ldots,x_{r}] $ be a monomial ideal, then it has unique minimal generating set, which we denote by $G(I)$.

\begin{theorem}The defining ideal $\mathfrak{P}_{nr}$ of monomial
curve associated to $\Gamma_{nr}$ is minimally generated by following binomials:
\medskip

\noindent $ f_{1}=x_{2}x_{3}^{3}-x_{1}x_{4}^{3}$;\\[2mm]
$ f_{(2,i)}=x_{1}^{n-i}x_{3}^{3i-1}-x_{2}^{n-i+1}x_{4}^{3i-2}, 1 \leqslant i \leqslant n$;\\[2mm]
$ f_{(3,j)}= x_{1}^{r-n+3+j}x_{2}^{n-1-j}-x_{3}^{2+3j}x_{4}^{r-1-3j}, 0 \leq j \leq n-1$;\\[2mm]
$ f_{(4,j)}=x_{1}^{r-2n+3+j}x_{2}^{2n-j}-x_{3}^{3j+1}x_{4}^{r+1-3j}, 0 \leqslant j \leqslant n-1$;\\[2mm]
$ f_{5}=x_{1}^{r-n+2}x_{2}^{n}x_{3}-x_{4}^{r+2}$;\\[2mm]
$ f_{6}=x_{2}^{n+1}x_{3}-x_{1}^{n}x_{4}^{2}$;\\[2mm]
$ f_{7}=x_{2}^{2n+1}-x_{1}^{2n-1}x_{3}x_{4}$.
\end{theorem}
 
\proof Let $J_{nr} = \langle\{f_{1}, f_{(2,i)}, f_{(3,j)},f_{(4,j)},f_{5},f_{6},f_{7}\mid 1\leq i\leq n, 0\leq j\leq n-1\}\rangle$, 
and  
\begin{align*}
\mathfrak{A}_{nr} = & 
\{ x_{1}x_{4}^{3}\} \cup \{x_{4}^{r+2}, x_{1}^{n}x_{4}^{2}, x_{2}^{2n+1}\} \cup \{x_{2}^{n-i+1}x_{4}^{3i-2} \mid 1\leq i\leq n\}\\ 
& \cup 
\{x_{1}^{r-n+3+j}x_{2}^{n-1-j}, x_{1}^{r-2n+3+j}x_{2}^{2n-j} \mid 0\leq j\leq n-1\} \cup \{x_{3}\}.
\end{align*}
Then, $J_{nr}+(x_{3})=\langle \mathfrak{A}_{nr}\rangle$ and it can be verified easily 
that $\mathfrak{A}_{nr}$ is minimal. First we note that 
$J_{nr}\subset \mathfrak{P}_{nr}$. We use the standard result: 
If $I=\langle G(I)\rangle$ is a monomial ideal in $A=k[x_{1},\ldots,x_{r}]$, then $A/I$ is a $k$-vector 
space whose basis consists of the images of monomials, which are not divisible by any elements 
of $G(I)$. Therefore $A/J_{nr}+\langle x_{3}\rangle $ is the vector space over $k$ and basis consists of 
the images of monomials, which are not divisible by any elements of $\mathfrak{A}_{nr}$. These are listed below:
\medskip

\noindent $S_{1}=\{x_{1}^{\alpha}:0 \leq \alpha \leq r+1\}$,\\[2mm]
$S_{2}=\{x_{2}^{\beta}:1 \leq \beta \leq 2n\}$,\\[2mm]
$S_{3}= \{x_{3}^{\gamma}:1 \leq \gamma \leq r+1\}$,\\[2mm]
$S_{4}=\{x_{1}^{\alpha}x_{2}^{\beta}:1 \leq \alpha \leq r-\beta+1,1\leq \beta \leq n-1 \}$,\\[2mm]
$S_{5}=\{x_{1}^{\alpha}x_{2}^{\beta}:1 \leq \alpha \leq r-\beta+2,n+1\leq \beta \leq 2n \}$,\\[2mm]
$S_{6}=\{x_{1}^{\alpha}x_{3}:1\leq \alpha \leq r+1 \}$,\\[2mm]
$S_{7}=\{x_{1}^{\alpha}x_{3}^{2}:1\leq \alpha \leq n-1 \}$,\\[2mm]
$S_{8}=\{x_{2}^{\beta}x_{3}^{\gamma} :1\leq \beta \leq n-1, 1 \leq \gamma \leq 3(n-\beta) \}$,\\[2mm]
$S_{9}=\{x_{1}^{\alpha}x_{2}^{\beta} x_{3}^{\gamma}:1\leq \alpha \leq n-1 , 1 \leq \beta \leq n-1 , 1 \leq \gamma \leq 2  \}$,\\[2mm]
$S_{10}=\{x_{1}^{\alpha}x_{\beta}x_{3}:n \leq \alpha \leq r-\beta +1 , 1 \leq \beta \leq n-1 \}$,\\[2mm]
$S_{11}=\{x_{1}^{\alpha}x_{2}^{n}:1\leq \alpha \leq r-n+2 \}$.
\medskip

It is clear from the expressions of the elements of $S_{i}$ that all the sets $S_{i}$, $1\leq i\leq 11$, 
are pairwise disjoint. Therefore the cardinality of this basis is,
\begin{align*}
&\displaystyle\sum_{i=1}^{11}\mid S_{i}\mid=(r+2)+2n+(r+1)+\frac{(n-1)}{2}(2r-n+2)+\frac{n}{2}(2r-3n+3)+(r+1)\\
&+(n-1)+\frac{3}{2}n(n-1)+2(n-1)^{2}+(r-n+2)+\frac{n-1}{2}(2r-3n+4)\\
& =3nr+3n+2r+4.
\end{align*}
Hence $\mathrm{dim}_{k} \left(A/J_{nr}+\langle x_{3}\rangle\right) =3nr+3n+2r+4=s+3n+1$ 
and by Theorem \ref{gastinger} it follows that $J_{nr}=\mathfrak{P}_{nr}$. 
\medskip

To show the minimality of the generating set, we consider the homomorphism 
\begin{align*}
\theta:k[x_{1},x_{2},x_{3},x_{4}]&\rightarrow k[x_{1},x_{2},x_{3},x_{4}]\\
\theta(x_{i})&=x_{i}, \, \mathrm{for}\,\, i=1,2,4;\\
\theta(x_{3})&=0.
\end{align*}
We note that, if $f$ is a generator of $J_{nr}$ or if $f=\sum_{g_{i} \in S \setminus \{f\}} p_{i}.g_{i}$, where $S$ is a 
generating set of $J_{nr}$, then $\theta(f)\in \mathfrak{A}_{nr}$. This implies that $\theta(f)=\sum_{ {g_{i} \in S \setminus \{f\}}} \theta(p_{i})\theta(g_{i})$, which gives a contradiction because $m \nmid m'$ for any pair of monomials 
$m,m' \in \mathfrak{A}_{nr}$. Therefore, $J_{nr}$ is a minimal generating set of $\mathfrak{P}_{nr}.$\qed 
\medskip

Let us denote the above generating set of $\mathfrak{P}_{nr}$ by 
$\mathfrak{S}_{nr}$, i.e.,
$$\mathfrak{S}_{nr}=\{f_{1}, f_{(2,i)}, f_{(3,j)},f_{(4,j)},f_{5},f_{6},f_{7}\mid 1\leq i\leq n, 0\leq j\leq n-1\}.$$

\begin{lemma}\label{gen2} 
Let $ g=x_{1}^{r+2}-x_{2}x_{4}^{r}$ be a polynomial in $k[x_{1},\ldots,x_{4}]$. Suppose $G_{nr}=(\mathfrak{S}_{nr}\setminus\{f_{(3,n-1)}\})\cup\{g\}$. Then $G_{nr}$ is also a generating set for the defining ideal $\mathfrak{P}_{nr}$.
\end{lemma}
\proof  Follows from the relation $ g=f_{(3,n-1)}+x_{4}^{r-3n+2}\cdot f_{(2,n)}$.\qed

\begin{theorem}\label{gbbc}
Let us consider the degree reverse lexicographic monomial order induced by $x_{1}>x_{2}>x_{3}> x_{4}$ 
on $k[x_{1},\ldots,x_{4}]$. Then, $G_{nr}$ is a Gr\"{o}bner basis of the defining ideal 
$\mathfrak{P}_{nr}$, 
with respect to the above order.
\end{theorem}

\proof We consider each $S$-polynomial and show that it reduces to zero 
upon division by $G_{nr}$. 

\begin{align*}
(1) \quad S(f_{1},f_{(2,1)})&=x_{2}^{n+1}x_{4}x_{3}-x_{1}^{n}x_{4}^{3}=x_{4}(x_{2}^{n+1}x_{3}-x_{1}^{n}x_{4}^{2})= x_{4}f_{6}.\\[2mm]
S(f_{1},f_{(2,i)})&=x_{2}^{n-i+2}x_{4}^{3i-2}- x_{1}^{n-i+1}x_{3}^{3i-4}x_{4}^{3}\\
&=-x_{4}^{3}(x_{1}^{n-i+1}x_{3}^{3i-4}-x_{2}^{n-i}x_{4}^{3i-5})
\\
&=-x_{4}^{3}f_{(2,i-1)},\quad \mathrm{for}\,\, 2\leq i\leq n.
\end{align*}

\begin{align*}
(2) \quad S(f_{1},f_{(3,n-2)})&=x_{3}^{3n-1}x_{4}^{r-3n+5}-x_{1}^{r+2}x_{4}^{3}\\
&=(x_{3}^{3n-1}-x_{2}x_{4}^{3n-2})x_{4}^{r-3n+5}-(x_{1}^{r+2}-x_{2}x_{4}^{r})x_{4}^{3}\\
&= f_{(2,n)}x_{4}^{r-3n+5}-gx_{4}^{3}.\\[2mm]
S(f_{1},f_{(3,i)})&=x_{3}^{5+3i}x_{4}^{r-1-3i}- x_{1}^{r-n+4+i}x_{2}^{n-2-i}x_{4}^{3}\\
&=-(x_{3}^{5+3i}x_{4}^{r-4-3i}- x_{1}^{r-n+4+i}x_{2}^{n-2-i})x_{4}^{3}\\
&=-f_{(3,i+1)}x_{4}^{3},\quad \mathrm{for}\,\, 2\leq i\leq n-3.
\end{align*}

\noindent $(3)$ \quad We consider two cases for $S(f_{1},f_{(4,i)})$:
\begin{enumerate}[(a)]
\item  For $i=n-1$,
\begin{align*}
S(f_{1},f_{(4,n-1)})&=x_{3}^{3n+1}x_{4}^{r+1-3n+3}-x_{1}^{r-n+3}x_{2}^{n}x_{4}^{3}
\\
&= (x_{3}^{3n-1}-x_{2}x_{4}^{3n-2})x_{3}^{2}x_{4}^{r-3n+4}-(x_{1}^{r-n+3}x_{2}^{n-1}-x_{3}^{2}x_{4}^{r-1})x_{2}x_{4}^{3}\\
&= f_{(2,n)}x_{3}^{2}x_{4}^{r-3n+4}-f_{(3,0)}x_{2}x_{4}^{3}.\\
\end{align*}

\item  For $ 0\leq i\leq n-2$,
\begin{align*}
S(f_{1},f_{(4,i)})&=x_{3}^{3i+4}x_{4}^{r+1-3i}- x_{1}^{r-2n+4+i}x_{2}^{2n-i-1}x_{4}^{3}\\
&=-(-x_{3}^{3i+4}x_{4}^{r-2-3i}+x_{1}^{r-2n+4+i}x_{2}^{2n-i-1})x_{4}^{3}=-f_{(4,i+1)}x_{4}^{3}.
\end{align*}
\end{enumerate}
\medskip

\noindent $(4) \quad S(f_{1},g)=x_{2}^{2}x_{3}^{3}x_{4}^{r}-x_{1}^{r+3}$. Since $\gcd(Lt(f_{1}),Lt(g))=1$, so $ S(f_{1},g) \longrightarrow_{G_{nr}} 0  $.
\medskip

\noindent $(5) \quad S(f_{1},f_{5})=-x_{1}^{r-n+3}x_{2}^{n-1}x_{4}^{3}+x_{4}^{r+2}x_{3}^{2}=-(x_{1}^{r-n+3}x_{2}^{n-1}-x_{3}^{2}x_{4}^{r-1})x_{4}^{3}=-f_{(3,0)}(x_{4}^{3})$.
\medskip

\noindent $(6) \quad S(f_{1},f_{6})=x_{1}^{n}x_{3}^{2}x_{4}^{2}-x_{1}x_{2}^{n}x_{4}^{3}
=(x_{1}^{n-1}x_{3}^{2}-x_{2}^{n}x_{4})x_{1}x_{4}^{2}=f_{(2,1)}(x_{1}x_{4}^{2})$.
\medskip

\noindent\begin{align*}
(7) \quad S(f_{1},f_{7}) & =  -x_{1}x_{2}^{2n}x_{4}^{3}+x_{1}^{2n-1}x_{3}^{4}x_{4}\\
{} & =  (x_{1}^{n-1}x_{3}^{2}-x_{2}^{n}x_{4})(x_{1}^{n}x_{3}^{2}x_{4}+x_{1}x_{4}^{2}x_{2}^{n})\\
{} & =  f_{(2,1)}(x_{1}^{n}x_{3}^{2}x_{4}+x_{1}x_{4}^{2}x_{2}^{n}).
\end{align*}

\noindent $(8)$ \quad 
For $ i<j$, we have 
\begin{align*}
S(f_{(2,i)},f_{(2,j)})&=x_{1}^{j-i}x_{2}^{n-j+1}x_{4}^{3j-2}-x_{2}^{n-i+1}x_{3}^{3(j-i)}x_{4}^{3i-2}\\
&=(x_{2}x_{3}^{3}-x_{1}x_{4}^{3})(\sum_{l=0}^{j-i-1}x_{1}^{l}x_{2}^{(n-i)-l}x_{3}^{3(j-i)-3(l+1)}x_{4}^{3i+(3l-2)} )\\
&=f_{1}(\sum_{l=0}^{j-i-1}x_{1}^{l}x_{2}^{(n-i)-l}x_{3}^{3(j-i)-3(l+1)}x_{4}^{3i+(3l-2)} ).
\end{align*}

\noindent $(9) \quad S(f_{(2,i)},f_{(3,j)})=x_{3}^{3(i+j)+1}x_{4}^{r-1-3j}-x_{1}^{r-2n+3+i+j}x_{2}^{2n-i-j}x_{4}^{3i-2}$. 
We consider three separate cases:
\begin{enumerate}[(a)]
\item For $i+j=n$,
\begin{align*}
 S(f_{(2,i)},f_{(3,j)})&=x_{3}^{3n+1}x_{4}^{r-1-3j}-x_{1}^{r-n+3}x_{2}^{n}x_{4}^{3i-2}\\
 &=(x_{1}^{r-n+3}x_{2}^{n-1}-x_{3}^{2}x_{4}^{r-1})(-x_{2}x_{4}^{3i-2})\\&+(x_{3}^{3n-1}-x_{2}x_{4}^{3n-2})(x_{3}^{2}x_{4}^{r-1-3j})
\\
&=f_{(3,0)}(-x_{2}x_{4}^{3i-2})+f_{(2,n)}(x_{3}^{2}x_{4}^{r-1-3j}).
\end{align*}

\item For $i+j>n$,
\begin{align*}
 S(f_{(2,i)},f_{(3,j)})&=f_{(3,(i+j)-n-1)}(-x_{1}x_{4}^{3i-2})+f_{(2,n)}(x_{3}^{3(i+j)-3n+2}x_{4}^{(r-1-3j)})\\&+f_{1}(x_{3}^{3(i+j)-3n-1}x_{4}^{r+3n-3j-3}).
 \end{align*}

\item For $i+j<n$,
\begin{align*}
 S(f_{(2,i)},f_{(3,j)})&=(x_{1}^{r-2n+3+i+j}x_{2}^{2n-i-j}-x_{3}^{3(i+j)+1}x_{4}^{r+1-3(i+j)})x_{4}^{3i-2}\\
 &=f_{(4,i+j)}x_{4}^{3i-2}.
\end{align*}
\end{enumerate}

\noindent $(10) \quad S(g,f_{(2,i)})=x_{1}^{r-n+2+i}x_{2}^{n-i+1}x_{4}^{3i-2}-x_{2}x_{3}^{3i-1}x_{4}^{r}=f_{(3,i-1)}x_{2}x_{4}^{3i-2}$,  \, for $1 \leq i \leq n-1$. For $i=n$, we note that 
$\gcd(Lt(g),Lt(f_{(2,n)})=1$, and hence $S(f_{(2,n)},g) \longrightarrow_{G_{nr}} 0$.
\medskip

\noindent $(11) S(f_{(2,i)},f_{(4,j)})=x_{3}^{3(i+j)}x_{4}^{r+1-3j}-x_{1}^{r-3n+3+i+j}x_{2}^{3n-i-j+1}x_{4}^{3i-2}$. 
We consider four separate cases:
\medskip
         
\begin{enumerate}[(a)]           
\item For $l=i+j \leq n-2$,
\begin{align*}
S(f_{(2,i)},f_{(4,j)})&=(x_{1}^{r-n+2+l}x_{2}^{n-l}-x_{3}^{3l-1}x_{4}^{r-3l+2})(-x_{3}x_{4}^{3i-1})+\\
&(x_{2}^{2n+1}-x_{1}^{2n-1}x_{3}x_{4})(-x_{1}^{r-3n+l+3}x_{2}^{n-l}x_{4}^{3i-2})\\
&=f_{(3,l-1)}(-x_{3}x_{4}^{3i-1})+f_{7}(-x_{1}^{r-3n+l+3}x_{2}^{n-l}x_{4}^{3i-2}).
\end{align*}

\item For $l=i+j=n-1$,
\begin{align*}
S(f_{(2,i)},f_{(4,j)})&=(x_{2}^{2n+1}-x_{1}^{2n-1}x_{3}x_{4}                                                                                                                                                      )(-x_{1}^{r-2n+2}x_{2}x_{4}^{3i-2})\\
&+(x_{1}^{r+1}x_{2}-x_{3}^{3n-4}x_{4}^{r-3n+5})(-x_{3}x_{4}^{3i-1})\\
&=f_{7}(-x_{1}^{r-2n+2}x_{2}x_{4}^{3i-2})+f_{(3,n-2)}(-x_{3}x_{4}^{3i-1}).\\
\end{align*}

\item For $l=i+j=n-1$,
\begin{align*}
S(f_{(2,i)},f_{(4,j)})&=(x_{1}^{r+2}-x_{2}x_{4}^{r})(-x_{3}x_{4}^{3i-1})\\
&+(x_{2}^{2n+1}-x_{1}^{2n-1}x_{3}x_{4}                                                                                                                             )(-x_{1}^{r-2n+3}x_{4}^{3i-2})\\&
+(x_{3}^{3n-1}-x_{2}x_{4}^{3n-2})(x_{3}x_{4}^{r+1-3j})\\
&=g(-x_{3}x_{4}^{3i-1})+f_{7}(-x_{1}^{r-2n+3}x_{4}^{3i-2})+f_{(2,n)}(x_{3}x_{4}^{r+1-3j}). \\
\end{align*}

\item For $ i+j\geq n+1$
\begin{align*}
S(f_{(2,i)},f_{(4,j)})&=f_{(4,(i+j)-n-1)}(-x_{1}x_{4}^{3i-2})+f_{(2,n)}(x_{3}^{3(i+j)-3n+1}x_{4}^{r+1-3j})\\
&+f_{1}(x_{3}^{3(i+j)-3n-2}x_{4}^{r-1-3j+3n}).
\end{align*} 
\end{enumerate}
 
\begin{align*}
(12) \quad S(f_{(2,i)},f_{5})=&x_{3}^{3i-2}x_{4}^{r+2}-x_{1}^{r-2n+2+i}x_{2}^{2n-i+1}x_{4}^{3i-2}\\
=&-(x_{1}^{r-2n+i+2}x_{2}^{2n-i+1}-x_{3}^{3i-2}x_{4}^{r+4-3i})x_{4}^{3i-2}\\
=& -f_{(4,i-1)}x_{4}^{3i+1},1 \leq i \leq n.
\end{align*}

\noindent $(13)$ \quad For $3 \leq i \leq n$,
\begin{align*}
S(f_{(2,i)},f_{6})&=x_{2}^{2n-i+2}x_{4}^{3i-2}-x_{1}^{2n-i}x_{4}^{2}x_{3}^{3i-2}\\
&=(x_{1}^{n-1}x_{3}^{2}-x_{2}^{n}x_{4})(-x_{1}^{n-i+1}x_{4}^{2}x_{3}^{3i-4})\\&-(x_{2}x_{3}^{3}-x_{1}x_{4}^{3})(\sum_{l=1}^{i-2}(x_{1}^{n-i+l}x_{2}^{(n-l)}x_{3}^{3i-(7+3(l-1))}x_{4}^{3l})\\
&=f_{(2,1)}(-x_{1}^{n-i+1}x_{4}^{2}x_{3}^{3i-4}-x_{2}^{n-i+2}x_{4}^{3i-3})\\&-f_{1}(\sum_{l=1}^{i-2}(x_{1}^{n-i+l}x_{2}^{(n-l)}x_{3}^{3i-(7+3(l-1))}x_{4}^{3l}).
\end{align*}

For $ i=2 $,
\begin{align*}S(f_{(2,2)},f_{6})&=x_{2}^{2n}x_{4}^{4}-x_{1}^{2n-2}x_{4}^{2}x_{3}^{4}\\
&=-(x_{1}^{n-1}x_{3}^{2}-x_{2}^{n}x_{4})(x_{1}^{n-1}x_{4}^{2}x_{3}^{2}+x_{2}^{n}x_{4}^{3})\\
&=-f_{(2,1)}(x_{1}^{n-i+1}x_{4}^{2}x_{3}^{3i-4}+x_{2}^{n-i+2}x_{4}^{3i-3}).
\end{align*}

For $ i=1$,  
\begin{align*}
S(f_{(2,1)},f_{6})=&(x_{2}^{2n+1}x_{4}-x_{1}^{2n-1}x_{4}^{2}x_{3})\\
=&(x_{2}^{2n+1}-x_{1}^{2n-1}x_{3}x_{4})x_{4}\\
=&f_{7}x_{4}.
\end{align*}
\medskip
  
\noindent $(14) \quad S(f_{(2,i)},f_{7})=-x_{1}^{n-i}x_{2}^{3n-i+2}x_{4}^{3i-2}+x_{1}^{3n-i-1}x_{4}x^{3i}$, \, 
 $1 \leqslant i \leqslant n$.\\[2mm]
 Since $\gcd(Lt(f_{(2,i)}),Lt(f_{7}))=1$, we have $ S(f_{(2,i)},f_{7}) \longrightarrow_{G_{nr}} 0  $.
\medskip

\begin{align*}
(15) \quad  S(f_{(3,i)},f_{(3,j)})&=-x_{3}^{2+3i}x_{4}^{r-1-3i}x_{1}^{j-i}+x_{2}^{j-i}x_{3}^{2+3j}x_{4}^{r-1-3j}\\
  &=(x_{2}x_{3}^{3}-x_{1}x_{4}^{3})(\sum_{l=0}^{j-i-1}x_{1}^{l}x_{2}^{(j-i)-(l+1)}x_{3}^{3j-(1+3l)}x_{4}^{r-3j+(-1+3l)})\\
 &=f_{1}(\sum_{l=0}^{j-i-1}x_{1}^{l}x_{2}^{(j-i)-(l+1)}x_{3}^{3j-(1+3l)}x_{4}^{r-3j+(-1+3l)}), \quad i<j.
\end{align*}

\allowdisplaybreaks

\begin{align*}
(16) \quad S(f_{(3,i),g})&=x_{2}^{n-i}x_{4}^{r}-x_{1}^{n-1-i}x_{3}^{2+3i}x_{4}^{r-1-3i}\\
&=-(x_{1}^{n-i-1}x_{3}^{3i+2}-x_{2}^{n-i}x_{4}^{3i+1})x_{4}^{r-1-3i}\\
&=- f_{(2,i+1)}x_{4}^{r-1-3i}. 
\end{align*}

\noindent $(17) \quad 
S(f_{(3,i)},f_{(4,j)})=x_{1}^{n+1-j}x_{3}^{3j+1}x_{4}^{r+1-3j}-x_{2}^{n-j+1+i}x_{3}^{2+3i}x_{4}^{r-1-3i}$. 
We consider three separate cases:
\medskip

\begin{enumerate}[(a)]
 \item For $i=j$,
\begin{align*}
S(f_{(3,i)},f_{(4,i)})&=(x_{2}x_{3}^{3}-x_{1}x_{4}^{3})(-x_{2}^{n}x_{3}^{3i-1}x_{4}^{r-1-3i}\\
&+(x_{1}^{n-1}x_{3}^{2}-x_{2}^{n}x_{4})x_{1}x_{3}^{3i-1}x_{4}^{r-3i+1}\\
&=f_{1}(-x_{2}^{n}x_{3}^{3i-1}x_{4}^{r-1-3i})+f_{(2,1)}x_{1}x_{3}^{3i-1}x_{4}^{r-3i+1}.
\end{align*}

\item For $i<j$, take $l=j-i$
\begin{align*}
S(f_{(3,i)},f_{(4,j)})&=(x_{1}^{n-l}x_{3}^{3l-1}-x_{2}^{n-l+1}x_{4}^{3l-2})x_{3}^{3i+2}x_{4}^{r+1-3j}\\
&=f_{(2,l)}x_{3}^{3i+2}x_{4}^{r+1-3j}.
\end{align*}

\item For $i>j$ and $l=j-i$, 
\begin{align*}
S(f_{(3,i)},f_{(4,j)})&=f_{(2,1)}x_{1}^{(i-j)+1}x_{3}^{3j-1}x_{4}^{r+1-3j}\\
&+(-f_{1}) \displaystyle\sum_{l=0}^{i-j}x_{1}^{l}x_{2}^{n-j+i-l}x_{3}^{3i-1-3l}x_{4}^{r-1-3i+3l}.
\end{align*}
\end{enumerate}

\begin{align*}
(18) \quad   S(f_{(3,i)},f_{5})=& x_{1}^{1+i}x_{4}^{r+2}-x_{2}^{1+i}x_{3}^{3+3i}x_{4}^{r-1-3i}\\
   &=-(x_{2}x_{3}^{3}-x_{1}x_{4}^{3})(\sum_{l=0}^{i}x_{1}^{l}x_{2}^{i-l}x_{3}^{3i-3l}x_{4}^{r-3i+3l-1})\\
   &=-f_{1}(\sum_{l=0}^{i}x_{1}^{l}x_{2}^{i-l}x_{3}^{3i-3l}x_{4}^{r-3i+3l-1}), 0 \leq i \leq n-2.
\end{align*}  
 
\begin{align*}
(19) \quad S( f_{(3,i)},f_{6})&=x_{1}^{r+3+i}x_{4}^{2}-x_{2}^{2+i}x_{3}^{3+3i}x_{4}^{r-1-3i}\\
=&(x_{1}^{r+2}-x_{2}x_{4}^{r})x_{1}^{1+i}x_{4}^{2}-(x_{2}x_{3}^{3}-x_{1}x_{4}^{3})(\sum_{l=0}^{i}x_{1}^{l}x_{2}^{1+(i-l)}x_{3}^{3(i-l)}x_{4}^{r-3i+3l-1})\\
&=g x_{1}^{1+i}x_{4}^{2}-f_{1}(\sum_{l=0}^{i}x_{1}^{l}x_{2}^{1+(i-l)}x_{3}^{3(i-l)}x_{4}^{r-3i+3l-1}), 0 \leq i \leq n-2.
\end{align*}
 
\begin{align*}
(20) \quad S(f_{(3,i)},f_{7})&= x_{1}^{r+n+2+i}x_{3}x_{4}-x_{2}^{n+i+2}x_{3}^{2+3i}x_{4}^{r-1-3i}\\
&=(x_{1}^{r+2}-x_{2}x_{4}^{r})(x_{1}^{n+i}x_{3}x_{4})+(x_{2}^{n+1}x_{3}-x_{1}^{n}x_{4}^{2})(-x_{1}^{i}x_{2}x_{3}x_{4}^{r-1})\\&-(x_{2}x_{3}^{3}-x_{1}x_{4}^{3})(\displaystyle\sum_{l=0}^{i-2}(x_{1}^{l}x_{2}^{n+i+(1-l)}x_{3}^{3i-(3l+1)}x_{4}^{r-3i-1+3l})\\
&=g(x_{1}^{n+i}x_{3}x_{4})+f_{6}(-x_{1}^{i}x_{2}x_{3}x_{4}^{r-1})\\
&-f_{1}(\displaystyle\sum_{l=0}^{i-2}(x_{1}^{l}x_{2}^{n+i+(1-l)}x_{3}^{3i-(3l+1)}x_{4}^{r-3i-1+3l}),\,\,0 \leq i < n-1.
\end{align*}

\begin{align*}
(21) \quad S(g,f_{(4,0)})&=-(x_{2}^{2n+1}-x_{1}^{2n-1}x_{3}x_{4})x_{4}^{r}\\
&=-f_{7}x_{4}^{r}.\\[2mm]
S(g,f_{(4,1)})&=(x_{1}^{n-1}x_{3}^{2}-x_{2}^{n}x_{4})(x_{1}^{n-i}x_{3}^{3i-1}x_{4}^{r+1-3i}+x_{2}x_{4}^{r-1})\\
&=f_{(2,1)}(x_{1}^{n-i}x_{3}^{3i-1}x_{4}^{r+1-3i}+x_{2}x_{4}^{r-1}).\\[2mm]
S(g,f_{(4,i)})&=x_{3}^{3i+1}x_{4}^{r+1-3i}-x_{1}^{j-2n+1}x_{2}^{2n-i-1}x_{4}^{r}\\
&=(x_{1}^{n-1}x_{3}^{2}-x_{2}^{n}x_{4})(x_{1}^{n-i}x_{3}^{3i-1}x_{4}^{r+1-3i}+x_{2}^{n-i+1}x_{4}^{r-1})\\
&+(x_{2}x_{3}^{3}-x_{1}x_{4}^{3})( \sum_{l=0}^{i-2}x_{1}^{n-i+l}x_{2}^{n-1-l}x_{3}^{3i-4-3l}x_{4}^{r+2-3i+3l})\\
&=f_{(2,1)}(x_{1}^{n-i}x_{3}^{3i-1}x_{4}^{r+1-3i}+x_{2}^{n-i+1}x_{4}^{r-1})\\
&+f_{1}( \sum_{l=0}^{i-2}x_{1}^{n-i+l}x_{2}^{n-1-l}x_{3}^{3i-4-3l}x_{4}^{r+2-3i+3l}),\quad i \geq 2.
\end{align*}
\medskip

\noindent $(22) \quad S(g,f_{5})=x_{1}^{n}x_{4}^{r+2}-x_{2}^{n+1}x_{3}x_{4}^{r}=(x_{2}^{n+1}x_{3}-x_{1}^{n}x_{4}^{2})(-x_{4}^{r})=f_{6}(-x_{4}^{r})$.
\medskip

\noindent $(23) \quad S(g,f_{6})=-x_{2}^{n+2}x_{3}x_{4}^{r}+x_{1}^{r+n+2}x_{4}^{2} \longrightarrow_{G} 0$, since 
$\gcd(Lt(g),Lt(f_{6}))=1$.
\medskip

\noindent $(24)$ \quad $ S(g,f_{7}) \longrightarrow_{G} 0$, since $\gcd(Lt(g),Lt(f_{7}))=1$. 
\medskip

\begin{align*}
(25) \quad  S(f_{(4,i)},f_{(4,j)})&=x_{2}^{j-i}x_{3}^{3j+1}x_{4}^{r+1-3j}-x_{1}^{j-i}x_{3}^{3i+1}x_{4}^{r+1-3i}\\
 &=(x_{2}x_{3}^{3}-x_{1}x_{4}^{3})(\sum_{l=0}^{j-i-1} x_{1}^lx_{2}^{(j-i)-(l+1)}x_{3}^{3j-(2+3l)}x_{4}^{r-3j+(1+3l)})\\
 &=f_{1}(\sum_{l=0}^{j-i-1} x_{1}^lx_{2}^{(j-i)-(l+1)}x_{3}^{3j-(2+3l)}x_{4}^{r-3j+(1+3l)}),\,\, i<j.
\end{align*}
\medskip

\begin{align*}
(26) \quad S(f_{(4,i)},f_{5})&=x_{2}^{n-i}x_{4}^{r+2}-x_{1}^{n-1-i}x_{3}^{3i+2}x_{4}^{r+1-3i}\\
 &=(x_{1}^{n-i-1}x_{3}^{3i+2}-x_{2}^{n-i}x_{4}^{3i-1})(-x_{4}^{r+1-3i})\\
 &=f_{(2,i+1)}(-x_{4}^{r+1-3i}),\,\,0 \leq i \leq n-1.
\end{align*}
\medskip

\noindent $(27)$ \quad We consider two cases for $S(f_{(4,i)},f_{6})$:

\begin{enumerate}[(a)]
\item For $0 \leqslant i < n-1$,
\begin{align*}S(f_{(4,i)},f_{6})&=x_{1}^{r-n+3+i}x_{2}^{n-i-1}x_{4}^{2}-x_{3}^{3i+2}x_{4}^{r+1-3i}\\
&= (x_{1}^{r-n+3+j}x_{2}^{n-1-j}-x_{3}^{2+3j}x_{4}^{r-1-3j})x_{4}^{2}
\\&= x_{4}^2f_{(3,i)}.
\end{align*}

\item For $i=n-1$, 
\begin{align*}
S(f_{(4,i)},f_{6})&=x_{1}^{r+2}x_{4}^{2}-x_{3}^{3n-1}x_{4}^{r-3n+4}\\&=(x_{1}^{r+2}-x_{2}x_{4}^{r})x_{4}^{2}-(x_{3}^{3n-1}-x_{2}x_{4}^{3n-2})x_{4}^{r-3n+4}\\
&=gx_{4}^{2}-f_{(2,n)}x_{4}^{r-3n+4}.
\end{align*}
\end{enumerate}

\noindent $(28)$ \quad $\mathrm{For}\,\, 0 \leqslant i \leqslant n-1,$\\ 
\begin{align*}
 S(f_{(4,i)},f_{7})&=x_{1}^{r+2+i}x_{3}x_{4}-x_{2}^{1+i}x_{3}^{3i+1}x_{4}^{r+1-3i}\\
 &=(x_{1}^{r+2}-x_{2}x_{4}^{r})(x_{1}^ix_{3}x_{4})\\
 &-(x_{2}x_{3}^{3}-x_{1}x_{4}^{3})(\sum_{l=0}^{i-1}x_{1}^{l}x_{2}^{i-l}x_{3}^{3i+(-3l-2)}x_{4}^{r-3i+1+3l})\\
 &= g(x_{1}^ix_{3}x_{4})-f_{1}(\sum_{l=0}^{i-1}x_{1}^{l}x_{2}^{i-l}x_{3}^{3i+(-3l-2)}x_{4}^{r-3i+1+3l}).
 \end{align*}
\medskip

\noindent $(29)$ \quad 
$S(f_{5},f_{6})=x_{1}^{r+2}x_{4}^{2}-x_{2}x_{4}^{r+2}=(x_{1}^{r+2}-x_{2}x_{4}^{r}) x_{4}^{2}=g x_{4}^{2}$.
\medskip

\begin{align*}
(30) \quad  S(f_{5},f_{7})=& -x_{2}^{n+1}x_{4}^{r+2}+x_{1}^{r+n+1}x_{3}^{2}x_{4}\\
 =&(x_{1}^{n-1}x_{3}^{2}-x_{2}^{n}x_{4})(x_{1}^{r+2}x_{4}+x_{2}x_{4}^{r+1})+( x_{1}^{r-n+3}x_{2}^{n-1}-x_{3}^{2}x_{4}^{r-1})x_{1}^{n-1}x_{2}x_{4}^{2}\\
 =& f_{(2,1)}(x_{1}^{r+2}x_{4}+x_{2}x_{4}^{r+1})+f_{(3,0)}x_{1}^{n-1}x_{2}x_{4}^{2}.
 \end{align*} 

\noindent $(31)$ \quad $S(f_{6},f_{7})= x_{1}^{2n-1}x_{3}^{2}x_{4}-x_{1}^{n}x_{2}^{n}x_{4}^{2}=(x_{1}^{n-1}x_{3}^{2}-x_{2}^{n}x_{4})x_{1}^{n}x_{4}= f_{(2,1)}x_{1}^{n}x_{4}$.
\medskip

\noindent Each $S$-polynomial reduces to zero, therefore, by Buchberger's Criterion,  
 $G_{nr}$ is a Gr\"{o}bner basis of $\mathfrak{P}_{nr}$, with respect to the degree 
 reverse lexicographic monomial order $>$ induced by $ x_{1}>x_{2}>x_{3}> x_{4}$. \qed 

\begin{corollary}\label{lead}
Let us consider the degree reverse lexicographic monomial order on $k[x_{1},\ldots,x_{4}]$, 
induced by $ x_{1}>x_{2}>x_{3}> x_{4}$. Then, with respect to this order, 
\begin{eqnarray*}
G(\mathrm{in}_{<}(\mathfrak{P}_{nr})) & = & \{x_{1}^{r+2}, x_{2}x_{3}^{3}, 
x_{1}^{r-n+2}x_{2}^{n}x_{3}, x_{2}^{n+1}x_{3}, x_{2}^{2n+1}\} \cup 
\{x_{1}^{n-i}x_{3}^{3i-1} \mid 1 \leq i \leq n\} \cup  \\ 
{} & {} & \{x_{1}^{r-n+3+j}x_{2}^{n-1-j} \mid 0 \leq j \leq n-2\} \cup \{x_{1}^{r-2n+3+l}x_{2}^{2n-l} 
\mid 0 \leq l \leq n-1\}.
\end{eqnarray*}
\end{corollary}
\proof Follows from Theorem \ref{gbbc}. \qed 
\medskip

We now use two theorems from \cite{hs}, written below, in order to determine the 
arithmetic Cohen-Macaulayness of the projective closure.
\medskip

\begin{lemma}\label{gbhom}
Let $I$ be an ideal in $A=k[x_{1},\ldots,x_{r}]$ and $I^{H}\subset A[x_{0}]$ its homogenization with respect to the variable $x_{0}$. Let $<$ be any reverse lexicographic monomial order on $A$ and $<_{0}$ the reverse lexicographic monomial order on $A[x_{0}]$ extended  from $A$ such that $x_{i}>x_{0}$  for all $i$.
\medskip

If $\{f_{1},\ldots,f_{n}\}$ is the reduced Gr\"{o}bner basis for $I$ w.r.t $<$, then $\{f_{1}^{H},\ldots,f_{n}^{H}\}$ is the reduced Gr\"{o}bner basis for $I^{H}$ w.r.t $<_{0}$, and $\mathrm{in}_{<_{0}}(I^{H})=(\mathrm{in}_{<}(I))A[x_{0}]$.
\end{lemma}

\proof See Lemma 2.1 in \cite{hs}.\qed

\begin{theorem}\label{criteria}
Let $\mathbf{n}:n_{1},\ldots,n_{r}$ be a sequence of positive integers with $n_{r}> n_{i}$ for all $i<n$. Let $<$ any reverse lexicographic order on $A=k[x_{1},\ldots,x_{r}]$ such that $x_{i}>x_{r}$ for all $1\leq i<r$. and $<_{0}$ the induced reverse lexicographic order on $A[x_{0}]$, where $x_{n}>x_{0}$. Then the following conditions are equivalent:
\begin{enumerate}
\item[(i)] The projective monomial curve $\overline{C(n_{1},\ldots,n_{r})}$ is arithmetically Cohen-Macaulay.
\item[(ii)] $\mathrm{in}_{<_{0}}((\mathfrak{p}(n_{1},\ldots,n_{r}))^{H})$ (homogenization with respect to $x_{0}$) is a Cohen-Macaulay ideal.
\item[(iii)] $\mathrm{in}_{<}(\mathfrak{p}(n_{1},\ldots,n_{r}))$ is a Cohen-Macaulay ideal.
\item[(iv)] $x_{r}$ does not divide any element of $G(\mathrm{in}_{<}(\mathfrak{p}(n_{1},\ldots,n_{r})))$.
\end{enumerate}
\end{theorem}

\proof See Theorem 2.2 in \cite{hs}. \qed

\begin{theorem}
The projective closure $\overline{B}_{nr}$ of the Backelin curve is arithmetically Cohen-Macaulay.
\end{theorem}

\proof From Corollary \ref{lead}, we see that $x_{4}$ does not divide any element of $G(\mathrm{in}_{<}(\mathfrak{P}_{nr})$. The proof follows from Theorem \ref{criteria}.\qed

\medskip
 
 \section{Hilbert Series of the Backelin curves}
 In this section, we compute the Hilbert series of the Backelin curve, using Algorithm 2.6 of \cite{coH}.
 \begin{lemma}\label{initial}
Let $I \subset k[x_1,...,x_n]$ be a graded ideal and $<$ a monomial order on
$k[x_1,...,x_n]$. Then $k[x_1,...,x_n]/I$ and $k[x_1,...,x_n]/ in_{<}(I)$ have the same Hilbert function, 
i.e. $$H(k[x_1,...,x_n]/I,i)= H(k[x_1,...,x_n]/in_{<}(I),i)$$ for all $i$.
\end{lemma}
\proof See Corollary 6.1.5 in \cite{hh}. \qed

\medskip
For a tuple $A=(a_{1},\ldots,a_{n})\in\mathbb{N}^{n}$, we write 
$x^{A} := x_{1}^{a_{1}}\cdots x_{n}^{a_{n}}$, and $|A|$ denotes the 
total degree of the monomial $x^{A}$.
\medskip

\begin{lemma}{\label{COH}}
Let $I=(x^{A_{1}},...,x^{A_{l}}) \subset k[x_1,...,x_n]$ be a monomial ideal. Let $p(I)$ denote the numerator of the Hilbert series of $ k[x_{1},x_{2},...,x_{n}]/I$. Then 
$$p(I)=p(x^{A_{1}})-\sum_{i=2}^{i=l}t^{|A_{i}|}p(x^{A_{1}},...,x^{A_{i-1}}:x^{A_{i}}),$$ 
where $p(x^{A_{1}},...,x^{A_{i-1}}:x^{A_{i}})$ 
denotes the numerator of the Hilbert Series of the ideal $(\langle x^{A_{1}},...,x^{A_{i-1}}\rangle :x^{A_{i}})$.
\end{lemma}

\proof See Corollary $2.3$ in \cite{coH}.\qed

\begin{theorem}
The numerator of the Hilbert series of the defining ideal $\mathfrak{P}_{nr}$ of 
the Backelin curve is 
\begin{align*}
&1-nt^{r+2}-2t^{r+3}+(3n+4)t^{r+4}-(2n+2)t^{r+5}-t^{2n+3}+2t^{2n+2}-t^{2n+1}+t^{n+4}+t^{n+3}\\&-t^{n+2}-t^{n+1}-t^{4}-\sum_{i=2}^{n}(t^{n+2i-1}+t^{n+2i+1}-2t^{n+2i}).
\end{align*}
\end{theorem}

\proof From Corollary \ref{lead}, we have,  $ \mathrm{in}_{<}(\mathfrak{P}_{nr})$ is generated by the set 
$$
\{x_{2}x_{3}^{3},x_{1}^{n-i}x_{3}^{3i-1},x_{1}^{r-n+3+j}x_{2}^{n-1-j},x_{1}^{r-2n+3+j} x_{2}^{2n-j},x_{1}^{r-n+2}x_{2}^{n}x_{3}, x_{2}^{n+1}x_{3},x_{2}^{2n+1}\mid 1 \leq i \leq n, 0 \leq j \leq n-1\}.
$$
We now use the algorithm 2.6 of \cite{coH}. Rearranging the generator of $\mathrm{in}_{<}(\mathfrak{P}_{nr})$, so that they are in ascending lexicographic 
order on the reversed set of variables $x_{4}>x_{3}>x_{2}>x_{1}$, we have 

\begin{align*}
G(\mathrm{in}_{<}(\mathfrak{P}_{nr}))=& \{x_{1}^{r+2},x_{1}^{r+1}x_{2},...,x_{1}^{r-n+3}x_{2}^{n-1},x_{1}^{r-n+2}x_{2}^{n+1},...,x_{1}^{r-2n+3}x_{2}^{2n},x_{2}^{2n+1},\\
& x_{1}^{r-n+2}x_{2}^{n}x_{3},x_{2}^{n+1}x_{3}, x_{1}^{n-1}x_{3}^{2},x_{1}x_{3}^{3},x_{1}^{n-2}x_{3}^{5},...,x_{3}^{3n-1}\}.
\end{align*}

\noindent By Lemma \ref{initial}, the Hilbert series of $k[x_{1},x_{2},x_{3},x_{4}]/\mathfrak{P_{nr}}$ 
is equal to the Hilbert series of $ \frac{k[x_{1},x_{2},x_{3},x_{4}]}{in_{<}(\mathfrak{P_{nr}})}$. 
Therefore, it is sufficient to compute the Hilbert series of the latter. 
Let $I$ denote the monomial ideal $in_{<}(\mathfrak{P_{nr}})$ and $p(I)$ is the numerator of the Hilbert series of  $in_{<}(\mathfrak{P_{nr}})$. Now using Lemma \ref{COH} to the ideal I, we have,
\begin{align*}
p(I)= & p(x_{1}^{r+2},x_{1}^{r+1}x_{2},...,x_{1}^{r-n+3}x_{2}^{n-1},x_{1}^{r-n+2}x_{2}^{n+1},...,x_{1}^{r-2n+3}x_{2}^{2n},x_{2}^{2n+1},\\
& x_{1}^{r-n+2}x_{2}^{n}x_{3},x_{2}^{n+1}x_{3},x_{1}^{n-1}x_{3}^{2},x_{2}x_{3}^{3},x_{1}^{n-2}x_{3}^{5},...,x_{3}^{3n-1}).
\end{align*}

\noindent Let \\[2mm]
$A_{1}=x_{1}^{r+2}$; $A_{2}=x_{2}^{2n+1}$; $A_{3}=x_{1}^{r-n+2}x_{2}^{n}x_{3}$; $A_{4}=x_{2}^{n+1}x_{3}$; $A_{5}=x_{2}x_{3}^{3}$;\\[2mm]
$B_{i}=x_{1}^{r-n+3+i}x_{2}^{n-1-i}$, $ 0 \leq i \leq n-2$;\\[2mm]
$C_{j}=x_{1}^{r-2n+3+j}x_{2}^{2n-j}$, $ 0 \leq j \leq n-1$;\\[2mm]
$D_{1}=x_{1}^{n-1}x_{3}^{2}$; $D_{l}=x_{1}^{n-l}x_{3}^{3l-1}$, $ 2 \leq l \leq n$. 
\medskip

\noindent After arranging them in the ascending lexicographic order on the reversed 
set of variables $x_{4}>x_{3}>x_{2}>x_{1}$, we have 
$$in_{<}(\mathfrak{P_{nr}})=\langle A_{1},B_{n-2},\dots,B_{0},C_{n-1},\dots,C_{0},A_{2},A_{3},A_{4},D_{1},A_{5},D_{2},\dots,D_{n}\rangle.$$
\allowdisplaybreaks

\noindent Using Lemma \ref{COH}, we now have,
\begin{align*}
p(I)&=p(A_{1})-t^{r+2}p(A_{1}:B_{n-2})-t^{r+2}\sum_{i=0}^{n-3} p(A_{1},B_{n-2},...,B_{i+1}:B_{i}) \\
&-t^{r+3}p(A_{1},B_{n-2},...,B_{0}:C_{n-1})-t^{r+3}\sum_{j=0}^{n-2} p(A_{1},B_{n-2},...,B_{0},C_{n-1},\dots,C_{j+1}:C_{j}) \\
&-t^{2n+1}p(A_{1},B_{n-2},...,B_{0},C_{n-1},\dots,C_{0}:A_{2})\\
&-t^{r+3}p(A_{1},B_{n-2},...,B_{0},C_{n-1},\dots,C_{0},A_{2}:A_{3})\\
&-t^{n+2}p(A_{1},B_{n-2},...,B_{0},C_{n-1},\dots,C_{0},A_{2},A_{3}:A_{4})\\
&-t^{n+1}p(A_{1},B_{n-2},...,B_{0},C_{n-1},\dots,C_{0},A_{2},A_{3},A_{4}:D_{1})\\
&-t^{4}p(A_{1},B_{n-2},...,B_{0},C_{n-1},\dots,C_{0},A_{2},A_{3},A_{4},D_{1}:A_{5})\\
&-\sum_{l=2}^{n}t^{n+2l-1}p(A_{1},B_{n-2},...,B_{0},C_{n-1},\dots,C_{0},A_{2},A_{3},A_{4},D_{1},A_{5},D_{2},\dots,D_{l-1}:D_{l});
\end{align*}

\noindent where
\begin{enumerate}[(i)]
\item $p(A_{1})=p(x_{1}^{r+2})=(1-t^{r+2})$;\\
\item $p(A_{1}:B_{n-2})=p(x_{1})=(1-t)$;\\
\item $p(A_{1},B_{n-2},...,B_{i+1}:B_{i})=p(x_{1})=(1-t)$;\\
\item $p(A_{1},B_{n-2},...,B_{0}:C_{n-1})=p(x_{1})=(1-t)$;\\
\item $p(A_{1},B_{n-2},...,B_{0},C_{n-1},\dots,C_{j+1}:C_{j})=p(x_{1})=(1-t)$;\\
\item  $p(A_{1},B_{n-2},...,B_{0},C_{n-1},\dots,C_{0}:A_{2})=p(x_{1}^{r-2n+3})=(1-t^{r-2n+3})$;\\
\item $p(A_{1},B_{n-2},...,B_{0},C_{n-1},\dots,C_{0},A_{2}:A_{3})=p(x_{1},x_{2})=(1+t^{2}-2t)$;\\
\item $p(A_{1},B_{n-2},...,B_{0},C_{n-1},\dots,C_{0},A_{2},A_{3}:A_{4})=p(x_{1}^{r-n+2},x_{1}^{r-2n+3}x_{2}^{n-1},x_{2}^{n})=(1-t^{n}-nt^{r-n+2}+nt^{r-n+3})$;\\
\item $p(A_{1},B_{n-2},...,B_{0},C_{n-1},\dots,C_{0},A_{2},A_{3},A_{4}:D_{1})=p(x_{2}^{n+1},x_{1}^{r-2n+3+j}x_{2}^{n-j})=(1-t^{n+1}-(n+1)t^{r-n+3}+(n+1)t^{r-n+4})$, where $ 0 \leq j \leq n$;\\
\item $p(A_{1},B_{n-2},...,B_{0},C_{n-1},\dots,C_{0},A_{2},A_{3},A_{4},D_{1}:A_{5})=p(x_{1}^{n-1},x_{2}^{n})=(1-t^{n-1}-t^{n}+t^{2n-1})$;\\[2mm]
\item $p(A_{1},B_{n-2},...,B_{0},C_{n-1},\dots,C_{0},A_{2},A_{3},A_{4},D_{1},A_{5},D_{2},\dots,D_{l-1}:D_{l})=p(x_{1},x_{2})=(1+t^{2}-2t)$.
\end{enumerate}

\noindent Therefore, 
\begin{align*}
p(I)&=(1-t^{r+2})-(n-1)t^{r+2}(1-t)-nt^{r+3}(1-t)-t^{2n+1}(1-t^{r-2n+3})-t^{r+3}(1+t^{2}-2t)\\&-t^{n+2}(1-t^{n}-nt^{r-n+2}+nt^{r-n+3})-t^{n+1}(1-t^{n+1}-(n+1)t^{r-n+3}+(n+1)t^{r-n+4})\\&-t^{4}(1-t^{n-1}-t^{n}+t^{2n-1})-(1+t^{2}-2t)(\sum_{i=2}^{n}t^{n+2i-1}),
\end{align*}
\noindent which is the same as
\begin{align*}
p(I)&=1-t^{r+2}-(n-1)t^{r+2}+(n-1)t^{r+3}-nt^{r+3}+nt^{r+4}-t^{2n+1}+t^{r+4}-t^{r+3}-t^{r+5}+2t^{r+4}\\&-t^{n+2}+t^{2n+2}+nt^{r+4}-nt^{r+5}-t^{n+1}+t^{2n+2}+(n+1)t^{r+4}-(n+1)t^{r+5}-t^{4}+t^{n+3}\\&+t^{n+4}-t^{2n+3}-\sum_{i=2}^{n}(t^{n+2i-1}+t^{n+2i+1}2t^{n+2i}).
\end{align*}
Therefore, 
\begin{align*}
p(I)&=1-nt^{r+2}-2t^{r+3}+(3n+4)t^{r+4}-(2n+2)t^{r+5}-t^{2n+3}+2t^{2n+2}-t^{2n+1}+t^{n+4}+t^{n+3}\\&-t^{n+2}-t^{n+1}-t^{4}-\sum_{i=2}^{n}(t^{n+2i-1}+t^{n+2i+1}-2t^{n+2i}).
\end{align*}
Hence, $\frac{p(I)}{(1-t)^4}$ is the Hilbert Series of the ring $k[x_{1},x_{2},x_{3},x_{4}]/I$.

\section{Syzygies of the affine and projective Backelin curve}

In this section, we compute the minimal free resolution of the defining ideal 
of the projective closure of the Backelin curve by computing the syzygies explicitly.

\begin{notations}\label{matrices B}{\rm 
Let us define the following matrices, which would be useful for writing the syzygies:
\begin{itemize}
\item[(a)]
\begin{align*}
\mathfrak{C^{1}_{H}}=( f_{1}^{H} , f_{(2,1)}^{H} , f_{6}^{H}, f_{(2,2)}^{H},...,f_{(2 ,\lceil \frac{n}{2})\rceil}^{H},f_{7}^{H},f_{(2,\lceil \frac{n}{2}\rceil+1)}^{H},...,f_{(2,n)}^{H},g^{H},f_{(3,n-2)}^{H},...,f_{(3,0)}^{H}\\,f_{(4,n-1)}^{H},...,f_{(4,0)}^{H}
,f_{5}^{H})
\end{align*}
\medskip

\item[(b)] 
$\mathfrak{C^{2}_{h}}:=
\collvec{&&R&&\\
0_{n \times n} & 0_{n \times (n+1)} & T^{i}_{n \times (n-2)} & T^{ii}_{n \times (n+1)} & T^{iii}_{n \times (n+2)}& T^{iv}_{n \times (n-1)}& T^{v}_{n \times 4}\\
S^{i}_{(n+2)\times n}& S^{ii}_{(n+2)\times(n+1)}& S^{iii}_{(n+2)\times (n-2)}&S^{iv}_{(n+2)\times (n+2)} &S^{v}_{(n+2)\times(n+2)}& S^{vi}_{(n+2)\times(n-1)}&S^{vii}_{(n+2)\times 4}\\
R^{i}_{(n+1) \times n} &R^{ii}_{(n+1) \times n}& 0_{(n+1) \times (n-2)} &R^{iii}_{(n+1) \times (n+1)} & R^{iv}_{(n+1) \times (n+2)} & R^{v}_{(n+1) \times (n-1)} & 0_{(n+1)\times 4}}_{(3n+4)\times (6n+5)}$
\medskip

\noindent such that,\\[2mm]
\noindent $\bullet R_{1 \times (6n+5)}=\begin{bmatrix}
r_{1,l}
\end{bmatrix}$, with \\[2mm]
$r_{1,1}=-x_{4}^{r-1}x_{0}$,\\[2mm]
$r_{1, 2+j}=-x_{3}^{1+3j}x_{4}^{3n-3j}x_{0}, \,\,\,0 \leq j \leq n-2$,\\[2mm]
$r_{1, n+2+i}=-x_{0}x_{3}^{2+3i}x_{4}^{r-3i-4}\,\,\,0 \leq i \leq n-3$,\\[2mm]
$r_{1, 2n+1+j}=-x_{2}^{j+1}x_{4}^{3n-5-3j}, \,\,\,0 \leq j \leq n-2$,\\[2mm]
$r_{1, 3n+2+j}=-x_{1}^{j}x_{3}^{3n-4-3j}, \,\,\,0 \leq j \leq n-2$,\\[2mm]
$r_{1, 4n}=-x_{1}^{n-2}x_{3}^{2}$,\\[2mm]
$r_{1, 6n+3}=-x_{1}^{n-1}$,\\[2mm]
$r_{1, 6n+5}=-x_{2}^{n}$,\\[2mm]
$r_{1,l}=0$, for \, $0 \leq j \leq n-2, \, 0 \leq i \leq n-3$ \, and\\[2mm]
$l \in \{1, \dots,6n+5\} \setminus 
\{1,2+j,n+2+i,2n+1+j,3n+2+j,4n,6n+3,6n+5\}$.\\[2mm]
\medskip

\noindent $\bullet T^{i}_{n \times (n-2)}=\begin{bmatrix}
0_{(\lceil \frac{n}{2}\rceil +2)\times (\lceil \frac{n}{2}\rceil -1)} & T^{i,1}_{(\lceil \frac{n}{2}\rceil +2)\times (\lceil \frac{n}{2}\rceil -1)}\\
T^{i,2}_{(\lceil \frac{n}{2}\rceil -2)\times (\lceil \frac{n}{2}\rceil -1)}&T^{i,3}_{(\lceil \frac{n}{2}\rceil -2)\times (\lceil \frac{n}{2}\rceil -1)}
\end{bmatrix}$
\medskip

\noindent $\bullet T^{i,1}_{(\lceil \frac{n}{2}\rceil +2)\times (\lceil \frac{n}{2}\rceil -1)}=\colvec{
0 & 0 & \cdots & 0  & 0 &-x_{3}^{3}   \\
0 & 0 & \cdots & 0  & 0 &0  \\
0 & 0 & \cdots & 0 & -x_{3}^{3}  & x_{1} \\
0 & 0 & \cdots & -x_{3}^{3}  & x_{1}& 0  \\
\vdots  & \vdots & \reflectbox{$\ddots$} & \reflectbox{$\ddots$} & \reflectbox{$\ddots$}& \vdots \\
0 &-x_{3}^{3}  & x_{1}& 0 & 0 & 0 \\
-x_{3}^{3} & x_{1} & 0 & 0 & 0 & 0 \\
 }$, 
\medskip
 
\noindent $\bullet T^{i,2}_{(\lceil \frac{n}{2}\rceil -2)\times (\lceil \frac{n}{2}\rceil -1)}=\colvec{
0 & 0 & \cdots & 0  & 0 &0 \\
0 & 0 & \cdots & 0  & 0 &-x_{3}^{3}   \\
0 & 0 & \cdots & 0 & -x_{3}^{3}  & x_{1} \\
0 & 0 & \cdots & -x_{3}^{3}  & x_{1}& 0  \\
\vdots  & \vdots & \reflectbox{$\ddots$} & \reflectbox{$\ddots$} & \reflectbox{$\ddots$}& \vdots \\
0 &-x_{3}^{3}  & x_{1}& 0 & 0 & 0 \\
-x_{3}^{3} & x_{1} & 0 & 0 & 0 & 0 \\
 }$
\medskip

\noindent $\bullet T^{i,3}_{(\lceil \frac{n}{2}\rceil -2)\times (\lceil \frac{n}{2}\rceil -1)}=\colvec{
0 &0 & 0 & \cdots &0 \\
x_{1} &0 & 0 & \cdots &0 \\
\vdots & \vdots & \vdots & \ddots &0\\
0 & 0& 0 & \cdots & 0\\
0 & 0 & 0 & \cdots & 0\\
}$
\medskip
 
\noindent $\bullet T^{ii}_{n \times (n+1)}=\begin{bmatrix}
T^{ii,1}_{(\lceil \frac{n}{2}\rceil +2)\times (\lceil \frac{n}{2}\rceil +1)} & T^{ii,2}_{(\lceil \frac{n}{2}\rceil +2)\times (\lceil \frac{n}{2}\rceil )}\\
T^{ii,3}_{(\lceil \frac{n}{2}\rceil -2)\times (\lceil \frac{n}{2}\rceil +1)}&T^{ii,4}_{(\lceil \frac{n}{2}\rceil -2)\times (\lceil \frac{n}{2}\rceil )}
\end{bmatrix}$
\medskip

\noindent $\bullet T^{ii,1}_{(\lceil \frac{n}{2}\rceil +2)\times (\lceil \frac{n}{2}\rceil +1)}=\colvec{
0 &0 & 0 & \cdots &0 \\
-x_{1}^{2n+4} &0 & 0 & \cdots &0 \\
\vdots & \vdots & \vdots & \ddots &0\\
0 & 0& 0 & \cdots & 0\\
0 & 0 & 0 & \cdots & 0\\
}$, 
\medskip

\noindent $\bullet T^{ii,2}_{(\lceil \frac{n}{2}\rceil +2)\times (\lceil \frac{n}{2}\rceil )}=\colvec{
0 & 0 & \cdots & 0  & 0 &-x_{4}^{3}   \\
0 & 0 & \cdots & 0  & 0 &0 \\
0 & 0 & \cdots & 0 & -x_{4}^{3} & x_{2} \\
0 & 0 & \cdots & -x_{4}^{3} & x_{2}& 0  \\
\vdots  & \vdots & \reflectbox{$\ddots$} & \reflectbox{$\ddots$} & \reflectbox{$\ddots$}& \vdots \\
0 &-x_{4}^{3} & x_{2}& 0 & 0 & 0 \\
-x_{4}^{3} & x_{2} & 0 & 0 & 0 & 0 \\
}$
\medskip

\noindent $\bullet T^{ii,3}_{(\lceil \frac{n}{2}\rceil -2)\times (\lceil \frac{n}{2}\rceil +1)}=\colvec{
0 & -x_{1}^{n+5}&0 & 0 & \cdots & 0  & 0 &0 \\
0 & 0&0 & 0 & \cdots & 0  & 0 &-x_{4}^{3}   \\
0 & 0&0 & 0 & \cdots & 0 & -x_{4}^{3}  & x_{2} \\
0 & 0&0 & 0 & \cdots &-x_{4}^{3}  & x_{2}& 0  \\
\vdots& \vdots& \vdots  & \vdots & \reflectbox{$\ddots$} & \reflectbox{$\ddots$} & \reflectbox{$\ddots$}& \vdots \\
0 & 0&0 &-x_{4}^{3}  & x_{2}& 0 & 0 & 0 \\
0 & 0&-x_{4}^{3} & x_{1} & 0 & 0 & 0 & 0 \\
 }$ 
\medskip

\noindent $\bullet T^{ii,4}_{(\lceil \frac{n}{2}\rceil -2)\times (\lceil \frac{n}{2}\rceil )}=\colvec{
0 &0 & 0 & \cdots &0 \\
x_{2} &0 & 0 & \cdots &0 \\ 
\vdots & \vdots & \vdots & \ddots &0\\
0 & 0& 0 & \cdots & 0\\
0 & 0 & 0 & \cdots & 0\\
}$
\medskip

\noindent $\bullet $ For $ 1 \leq i \leq (n-1)$,\\[2mm]
$T^{iii}_{n \times (n+2)}=
\colvec{&1&  &i^{th}& & n-1\\
 -x_{1}^{r-2n+3}x_{2}^{n} & 0 &\cdots &0 &\cdots  & 0 & 0 & x_{1}^{n}x_{4}\\
0 & -x_{1}^{r-2n+3}x_{2}^{n-1}  &\cdots &-x_{1}^{r-2n+2+i}x_{2}^{n-i}& \cdots & -x_{1}^{r-n+1}x_{2} &-x_{1}^{r-n+2}& -x_{2}^{n}\\

      0 & 0 &\cdots &0& \cdots & 0& 0 & x_{3}& (\lceil \frac{n}{2}\rceil +3)^{th} \\
  0 & 0 &\cdots &0& \cdots & 0& 0 & 0 & row \\
  \vdots & \vdots &\ddots &\vdots& \ddots & \vdots& \vdots & \vdots \\
   0 & 0 &\cdots &0& \cdots & 0& 0 & 0 \\     
}$
\medskip

\noindent $\bullet $ For $ 0 \leq i \leq (n-2)$,\\[2mm]
$T^{iv}_{n \times (n-1)}=
\colvec{&  &i^{th}&\\
 x_{1}^{2n+4}x_{2}  &\cdots &x_{1}^{2n+4-i}x_{2}^{i+1}& \cdots & x_{1}^{n+6}x_{2}^{n-1} \\

  0 &\cdots &0& \cdots & 0 \\
   0 &\cdots &0& \cdots & 0 \\
  \vdots &\ddots &\vdots& \ddots  \\
  0 &\cdots &0& \cdots & 0  \\     
}$
\medskip

\noindent $\bullet T^{v}_{n \times 4}=
\collvec{
x_{2}^{n+1} & -x_{2}x_{3} & -x_{1}^{r-n+3}& x_{1}x_{4}^{2}\\
-x_{1}^{n-1}x_{3}& x_{4} & 0 & x_{3}^{2}\\
0 & 0 & 0 &0\\
\vdots & \vdots & \vdots & \vdots\\
x_{4} & 0 & 0 &0 & (\lceil \frac{n}{2}\rceil +3)^{th}\\
0 & 0 & 0 &0\\
\vdots & \vdots & \vdots & \vdots\\
0 & 0 & 0 &0 
}$
\medskip

\noindent $\bullet S^{i}_{(n+2)\times n}=\collvec{
0 &0 & 0 & \cdots &0 \\
\vdots & \vdots & \vdots & \ddots &0\\
0 & 0& 0 & \cdots & 0\\
-x_{2}x_{3} & 0 & 0 & \cdots & 0\\
}$
\medskip

\noindent $\bullet S^{ii}_{(n+2)\times(n+1)}=\collvec{
x_{3}^{2}x_{4}^{r-3n+4}x_{0}&0 & 0 & \cdots & 0 &  x_{0}x_{4}^{r-3n+5} & -x_{3}^{3} \\
0 &0& 0 & \cdots & 0 &  0 & x_{1}\\
0 & 0&0 & \cdots & 0 &  -x_{2} & 0\\
0 & 0&0 &\cdots & -x_{2}& x_{1}&0 \\
0 & 0&0 &\reflectbox{$\ddots$} & x_{1}& 0&0 \\
\vdots &\vdots  & \reflectbox{$\ddots$} & \reflectbox{$\ddots$} & \reflectbox{$\ddots$}& \vdots & \vdots\\
0 &-x_{2}&  x_{1} & 0& 0 & 0 & 0 \\
 -x_{2}^{2}& x_{1} & 0 & 0 & 0 & 0
}$
\medskip

\noindent $\bullet S^{iii}_{(n+2)\times (n+1)}=\collvec{
x_{1} &0 & 0 & \cdots &0 \\
\vdots & \vdots & \vdots & \ddots &0\\
0 & 0& 0 & \cdots & 0\\
0 & 0 & 0 & \cdots & 0\\
}$
\medskip

\noindent $\bullet S^{iv}_{(n+2)\times (n+1)}=\collvec{
0 & 0 & -x_{4}^{3} & x_{2} & 0 & \cdots & 0\\
0 & 0 & x_{2} & 0& 0 & \cdots & 0\\
-x_{4}^{2} & -x_{3}x_{4}& 0  & 0& 0 & \cdots & 0\\
0 &0 & 0 & 0 & 0& \cdots &0 \\
\vdots & \vdots &\vdots& \vdots&\vdots & \ddots &0 \\
0 &0 & 0 & 0 & 0& \cdots &0 \\
}$
\medskip

\noindent $\bullet S^{v}_{(n+2)\times(n+1)}=\collvec{
0 & 0 & \cdots & 0 &  0\\
0 & 0 & \cdots & 0 &0 \\ 
0 & 0 & \cdots & 0 &  -x_{4}^{2}\\
0 & 0 & \cdots &   -x_{4}^{2} & 0  \\
\vdots  & \vdots & \reflectbox{$\ddots$} & \reflectbox{$\ddots$} & \reflectbox{$\ddots$} \\
0 &  -x_{4}^{2} & 0& 0 & 0 \\
}$
\medskip

\noindent $\bullet S^{vi}_{(n+2)\times(n-1)}=\collvec{
0 & 0 & \cdots & 0 &  0\\
0 & 0 & \cdots & 0 &0  \\ 
0 & 0 & \cdots & 0 &  -x_{3}^{2} \\
0 & 0 & \cdots &   -x_{3}^{2} & 0  \\
\vdots  & \vdots & \reflectbox{$\ddots$} & \reflectbox{$\ddots$} &\vdots \\
0 &  -x_{3}^{2} & 0& 0 & 0  \\
 -x_{3}^{2}& 0 & 0& 0 & 0  \\
 0 &  0 & 0& 0 & 0  \\
}$
\medskip

\noindent $\bullet S^{vii}_{(n+2)\times 4}= \collvec{
0 & 0 & 0 &0 \\
0 &0 &0 &0\\
0 & 0 & x_{3}^{2} & 0\\
0 & 0 & 0 & 0\\
\vdots & \vdots & \vdots & \vdots\\
0 & 0& 0  & 0\\
0 & 0 & -x_{2}x_{4}& 0\\
}$ 
\medskip

\noindent $\bullet R^{i}_{(n+1)\times n}=\collvec{
0 & 0 & \cdots & 0  & 0 &-x_{2}   \\
0 & 0 & \cdots & 0 & -x_{2}  & x_{1} \\
0 & 0 & \cdots & -x_{2}  & x_{1}& 0  \\
\vdots  & \vdots & \reflectbox{$\ddots$} & \reflectbox{$\ddots$} & \reflectbox{$\ddots$}& \vdots \\
0 &-x_{2}  & x_{1}& 0 & 0 & 0 \\
0 & x_{1} & 0 & 0 & 0 & 0 \\
 x_{1}& 0 & 0 & 0 & 0 & 0
}$
\medskip

\noindent $\bullet R^{ii}_{(n+1) \times (n+1)}=\collvec{
x_{1} & \cdots & 0 & 0 &0\\
0& \cdots & 0 & 0 &0\\
\vdots & \vdots &\vdots &\vdots &\vdots \\
0 & \cdots & 0 & 0 &0\\
0 & \cdots & 0 &0 &0\\
}$
\medskip

\noindent $\bullet R^{iii}_{(n+1)\times (n+1)}=\collvec{
0 &0 & 0 & \cdots &0 \\
\vdots & \vdots & \vdots & \ddots &0\\
0 &0 & 0 & \cdots &0 \\
0 & x_{2}& 0 & \cdots & 0\\
x_{2} & 0 & 0 & \cdots & 0\\
}$
\medskip

\noindent $\bullet R^{iv}_{(n+1)\times (n+2)}=\collvec{
0 & 0 & \cdots & 0 &  x_{3} & 0 \\
0 & 0 & \cdots &  x_{3} & 0& 0  \\
\vdots  & \vdots & \reflectbox{$\ddots$} & \reflectbox{$\ddots$} & \reflectbox{$\ddots$}& \vdots \\
-x_{4} & x_{3} & 0& 0 & 0 & 0 \\
 x_{3}& 0 & 0 & 0 & 0 & 0
}$
\medskip

\noindent $\bullet R^{v}_{(n+1)\times(n-1)}=\collvec{
0 & 0 & \cdots & 0 & x_{4} \\
0 & 0 & \cdots & x_{4} & 0  \\
\vdots  & \vdots & \reflectbox{$\ddots$} & \reflectbox{$\ddots$} & \vdots \\
0 & x_{4} & 0& 0 & 0 \\
x_{4}& 0 & 0 & 0 & 0 \\
0& 0 & 0 & 0 & 0
}$
\medskip

\item[(c)]
$\mathfrak{C^{3}_{h}}:= \mathfrak{C^{3}}_{(6n+5)\times (3n+2)}=\begin{bmatrix}
X_{(n+1)\times (n+1)} & Y_{(n+1)\times (n-1)} &  & 0_{(n+1) \times (n+2)}\\
X'_{(n-1)\times (n+1)} & Y'_{(n-1)\times (n-1)} &  &0_{(n-1) \times (n+2)}\\
X''_{n\times (n+1)} & 0_{(n+1) \times (n-1)} & &Z''_{n\times (n+2)}\\
X^{iii}_{n\times (n+1)} & Y^{iii}_{n\times (n-1)} & &Z^{iii}_{n\times (n+2)}\\
X^{iv}_{(n+1)\times (n+1)} & Y^{iv}_{(n+1)\times (n-1)} & &Z^{iv}_{(n+1)\times (n+2)}\\
X^{v}_{(n+1)\times (n+1)} & Y^{v}_{(n+1)\times (n-1)} & & 0_{(n+1)\times (n+2)}\\
X^{vi}_{4\times (n+1)} & 0_{4\times (n-1)} & &Z^{vi}_{4\times (n+2)}\\
\end{bmatrix}
$,\\[2mm]
where\\[2mm]
\noindent $\bullet X_{(n+1)\times (n+1)}=\collvec{
0 & 0 & 0 &\cdots &0 &x_{2} & 0\\
0 & 0 & 0  &\cdots & 0 & 0 & 0\\
0 & 0 & 0 & \cdots & -x_{4} & 0 & 0\\
\vdots & \vdots & \reflectbox{$\ddots$} &  \reflectbox{$\ddots$} &  \reflectbox{$\ddots$} & \vdots & \vdots\\
0 & 0 & -x_{4} & \reflectbox{$\ddots$} & 0 & 0 & 0\\
0 & -x_{4} & 0 & \hdots & 0 & 0 &-x_{3}\\
-x_{4} & 0 & 0 & \hdots  & 0 & -x_{3} & 0 \\
}
$

\noindent $\bullet Y_{(n+1)\times (n-1)}= \collvec{
0 & 0 & \cdots & 0  & 0 &-x_{3} \\
0 & 0 & \cdots & 0 & x_{3} & x_{4}\\
0 & 0 & \cdots & x_{3} & 0 & 0\\
\vdots  & \vdots & \reflectbox{$\ddots$} & \reflectbox{$\ddots$} & \reflectbox{$\ddots$}& \vdots\\
0 & x_{3} & 0 & 0 & 0 & 0\\
x_{3} & 0 & 0 & 0 & 0 & 0\\
0 & 0 & 0 & 0 & 0 & 0\\
0 & 0 & 0 & 0 & 0 & 0\\

}
$

\noindent $\bullet X'_{(n-1)\times (n+1)}=\collvec{
0 & 0 & 0 & \cdots & x_{3}^{2} & 0 & 0\\
\vdots & \vdots & \reflectbox{$\ddots$} &  \reflectbox{$\ddots$} &  \reflectbox{$\ddots$} & \vdots & \vdots\\
0 & 0 & x_{3}^{2} & \reflectbox{$\ddots$} & 0 & 0 & 0\\
0 & x_{3}^{2} & 0 & \hdots & 0 & 0 &0\\
x_{3}^{2} & 0 & 0 & \hdots  & 0 & 0 & x_{4}^{2}\\
}$
\medskip

\noindent $\bullet Y'_{(n-1)\times (n+1)}=\collvec{
0 & 0 & 0 & \cdots & 0 &x_{4}^{2}\\
0 & 0 & 0 & \cdots & x_{4}^{2}& 0\\
\vdots & \vdots & \reflectbox{$\ddots$} &  \reflectbox{$\ddots$} &  \reflectbox{$\ddots$} & \vdots\\
0 & 0 & x_{4}^{2} & \reflectbox{$\ddots$} & 0 & 0\\
0 & x_{4}^{2} & 0 & \hdots & 0 & 0\\
x_{4}^{2} & 0 & 0 & \hdots  & 0 & 0 \\
}$

\noindent $\bullet X''_{n \times (n+1)}=\collvec{
0 & \cdots & 0 & 0 &0\\
\vdots & \vdots &\vdots &\vdots &\vdots \\
0 & \cdots & 0 & 0 &0\\
0 & \cdots & 0 &-x_{1} &0\\
0 & \cdots & 0 & 0 &0
}$
\medskip

\noindent $\bullet Z''_{n \times (n+2)}=\collvec{
0 & 0 & \cdots & 0  & 0 &-x_{2} & 0 & 0 &0 \\
0 & 0 & \cdots & 0 & -x_{2} & x_{4}^{3}& 0 & 0 &0\\
0 & 0 & \cdots &-x_{2} & x_{4}^{3} & 0 & 0 & 0 &0\\
\vdots  & \vdots & \reflectbox{$\ddots$} & \reflectbox{$\ddots$} & \reflectbox{$\ddots$}& \vdots &\vdots& \vdots& \vdots\\
0 &-x_{2} & x_{4}^{3} & 0 & 0 & 0 & 0 & 0 & 0\\
-x_{2} & x_{4}^{3} & 0 & 0 & 0 & 0 & 0 & 0 & 0\\
0 & 0 & 0 & 0 & 0 & 0 & 0 &-x_{3} & 0\\
0 & 0 & 0 & 0 & 0 & 0 & x_{1} & x_{4} & 0\\
}$
\medskip

\noindent $\bullet X^{iii}_{n \times (n+1)}=\collvec{
0 & \cdots & 0 & 0 &x_{4}^{r-3n+4}x_{0}\\
0 & \cdots & 0 & 0 &0\\
\vdots & \vdots &\vdots &\vdots &\vdots \\
0 & \cdots & 0 & 0 &0\\
0 & \cdots & 0 &-x_{4}^{r-3n+4}x_{0} &0\\
}$
\medskip

\noindent $\bullet Y^{iii}_{n \times (n+1)}=\collvec{
0 & \cdots & 0 & 0 &0\\
\vdots & \vdots &\vdots &\vdots &\vdots \\
0 & \cdots & 0 & 0 &0\\
0 & \cdots & 0 &0 &x_{1}\\
}$
\medskip

\noindent $\bullet Z^{iii}_{n \times (n+2)}=\collvec{
0 & 0 & \cdots & 0  & 0 &x_{1} & 0 & 0 &0 \\
0 & 0 & \cdots & 0 & x_{1} & -x_{3}^{3}& 0 & 0 &0\\
0 & 0 & \cdots & x_{1} & -x_{3}^{3}& 0 & 0 & 0 &0\\
\vdots  & \vdots & \reflectbox{$\ddots$} & \reflectbox{$\ddots$} & \reflectbox{$\ddots$}& \vdots &\vdots& \vdots& \vdots\\
0 & x_{1} & -x_{3}^{3}& 0 & 0 & 0 & 0 & 0 & 0\\
 x_{1} & -x_{3}^{3}& 0 & 0 & 0 & 0 & 0 & 0 & 0\\
0 & 0 & 0 & 0 & 0 & 0 & 0 &x_{2} & 0
}$
\medskip

\noindent $\bullet X^{iv}_{(n+1) \times (n+1)}=\collvec{
0 & \cdots & 0 & 0 &0\\
\vdots & \vdots &\vdots &\vdots &\vdots \\
0 & \cdots & 0 & 0 &0\\
0 & \cdots & 0 &0 &x_{1}\\
0 & \cdots & 0 &x_{1} & -x_{2}\\
0 & \cdots & 0 & 0 &0\\
}$
\medskip

\noindent $\bullet Y^{iv}_{(n+1) \times (n+1)}=\collvec{
0 & 0 & \cdots & 0  & 0 &x_{1} & 0 \\
0 & 0 & \cdots & 0 & x_{1} & -x_{2}& 0 \\
0 & 0 & \cdots &x_{1} & -x_{2} & 0 & 0 \\
\vdots  & \vdots & \reflectbox{$\ddots$} & \reflectbox{$\ddots$} & \reflectbox{$\ddots$}& \vdots &\vdots\\
0 &x_{1} & -x_{2} & 0 & 0 & 0 & 0 \\
x_{1} & -x_{2} & 0 & 0 & 0 & 0 & 0 \\
-x_{2} & 0 & 0 & 0 & 0 & 0 & 0 \\
0 & 0 & 0 & 0 & 0 & 0 & 0 \\
0 & 0 & 0 & 0 & 0 & 0 & 0 \\
}$
\medskip

\noindent $\bullet Z^{iv}_{(n+1) \times (n+2)}=\collvec{
0 & \cdots & 0 & 0 &x_{2} & 0  & -x_{4}\\
0 & \cdots & 0 & 0 &0 & 0 & 0\\
\vdots & \vdots &\vdots &\vdots &\vdots &\vdots &\vdots\\
0 & \cdots & 0 & 0 &0 & 0 & 0\\
0 & \cdots & 0 &0 &-x_{1}^{r-2n+3} & 0 &0\\
}$
\medskip

\noindent $\bullet X^{v}_{(n-1) \times (n+1)}=\collvec{
0 & 0 & \cdots & 0  & 0 &x_{1} & 0 & 0  \\
0 & 0 & \cdots & 0 & x_{1} & -x_{2}& 0 & 0 \\
0 & 0 & \cdots & x_{1} & -x_{2}& 0 & 0 & 0  \\
\vdots  & \vdots & \reflectbox{$\ddots$} & \reflectbox{$\ddots$} & \reflectbox{$\ddots$}& \vdots &\vdots & \vdots\\
0 & x_{1} & -x_{2}& 0 & 0 & 0 & 0 & 0\\
 x_{1} & -x_{2}& 0 & 0 & 0 & 0 & 0 & 0\\
}$
\medskip

\noindent $\bullet Y^{v}_{(n-1) \times (n+1)}=\collvec{
0 & \cdots & 0 & 0 &x_{2}\\
\vdots & \vdots &\vdots &\vdots &\vdots \\
0 & \cdots & 0 & 0 &0\\
0 & \cdots & 0 &0 &0\\
}$
\medskip

\noindent $\bullet X^{vi}_{4 \times (n+1)}=\collvec{
0 & 0 & \cdots 0 \\
0 & 0 & \cdots 0 \\
x_{2} & 0 & \cdots 0 \\
0 & 0 & \cdots 0 \\
}$
\medskip

\noindent $\bullet Z^{vi}_{4 \times (n+2)}=\collvec{
0 & 0 & \cdots &0 & x_{1}^{r-2n+3} & x_{3}\\
-x_{3}^{2} & 0 & \cdots &0 & 0 & -x_{2}^{n}\\
0 & 0 & \cdots &-x_{4} & 0 & 0\\
x_{4} & 0 & \cdots &0& 0 & x_{1}^{n-1}\\
}$.
\end{itemize}
}
\end{notations}
\medskip

Let us recall the Buchsbaum Eisenbud acyclicity criterion.

\begin{lemma}\label{Exact}
Let $R$ be a Noetherian ring, and 
\[
F.\,: 0 \rightarrow F_{s} \xrightarrow{\phi_{s}} F_{s-1} \xrightarrow{} \cdots \xrightarrow{}F_{1} \xrightarrow{\phi_{1}} F_{0} \rightarrow 0,
\]
a complex of finite free $R$-modules. Set $r_{i}=\sum_{j=i}^{s}(-1)^{j-i}\mathrm{rank} F_{j}$. 
Then the following statements are equivalent:
\begin{enumerate}[(a)]
\item $F.$ is acyclic;
\item $\mathrm{grade}\,I_{r_{i}}(\phi_{i})\geq i$ for $ i=1,\dots,s$.
\end{enumerate}
\proof See Theorem 1.4.13 in \cite{acyclicity}.

\end{lemma}

\medskip

\begin{notations}{\rm 
Given an $n \times m$ matrix $A$ with entries $a_{ij}$, a minor of $A$ is denoted by
$$M=[r_{1} \,\,r_{2}\,\,\cdots\,\,r_{p} \vert c_{1} \,\,c_{2}\,\,\cdots\,\,c_{p}]=
\mathrm{det}\begin{bmatrix}
a_{r_{1}c_{1}} & a_{r_{1}c_{2}}& \cdots &a_{r_{1}c_{p}}\\
\vdots & \vdots &\vdots &\vdots \\
a_{r_{p}x_{1}} & a_{r_{p}c_{2}}& \cdots &a_{r_{p}c_{p}}
\end{bmatrix},$$ 
where $\{ r_{1},\dots.r_{p}\} \subseteq \{1, \dots,n\}$ and 
$\{ c_{1},\dots,c_{p}\} \subseteq \{1, \dots,m\}$. 
}
\end{notations}

\begin{theorem}
For $ n \geq 4$, a minimal graded free resolution of the defining ideal of the projective 
closure of the Backelin curve  $\overline{\Gamma}_{nr}$\, is given by
$$ \overline{\mathfrak{R}_{h}}: 0 \rightarrow R^{3n+2} \xrightarrow{\mathfrak{C^{3}_{h}}} R^{6n+5} \xrightarrow{\mathfrak{C^{2}_{h}}} R^{3n+4} \xrightarrow{\mathfrak{C^{1}_{h}}} R \rightarrow R/\overline{\mathfrak{P_{nr}}}\rightarrow 0,$$
where the matrices $\mathfrak{C^{i}_{h}}, 1 \leq i \leq 3 $ are defined in notation \ref{matrices B}. 
Moreover, for $ n \geq 4$, a minimal graded free resolution of the defining 
ideal of the Backelin curve  $\Gamma_{nr}$ is given by,
$$ \mathfrak{R}_{h}: 0 \rightarrow R^{3n+2} \xrightarrow{\mathfrak{\tilde{C}^{3}_{h}}} R^{6n+5} \xrightarrow{\mathfrak{\tilde{C}^{2}_{h}}} R^{3n+4} \xrightarrow{\mathfrak{\tilde{C}^{1}_{h}}} R \rightarrow R/\mathfrak{P_{nr}}\rightarrow 0,$$
where the matrices $\mathfrak{\tilde{C}^{i}_{h}}, 1 \leq i \leq 3 $, are obtained by evaluating $x_{0}=1$ in  $\mathfrak{C^{i}_{h}}, 1 \leq i \leq 3 $.
\end{theorem}

\proof  We can check that  $$ \overline{\mathfrak{R}_{h}}: 0 \rightarrow R^{3n+2} \xrightarrow{\mathfrak{C^{3}_{h}}} R^{6n+5} \xrightarrow{\mathfrak{C^{2}_{h}}} R^{3n+4} \xrightarrow{\mathfrak{C^{1}_{h}}} R \rightarrow R/\overline{\mathfrak{P_{nr}}}\rightarrow 0$$
is a chain complex by calculating $\mathfrak{C^{i}_{h}}\circ \mathfrak{C^{i+1}_{h}}=0$, $1 \leq i \leq 2$. 
To prove the exactness of the above chain complex, we use the Buchsbaum-Eisenbud acyclicity criterion (see Lemma \ref{Exact}).
Let $ r_{i}$ be the $i^{th}$ expected rank of $\mathfrak{C^{i}_{h}}$. 
Then, 
\begin{itemize}
\item $r_{1}=\mathrm{rank}(\mathfrak{C^{1}_{h}})-\mathrm{rank}(\mathfrak{C^{2}_{h}})+\mathrm{rank}(\mathfrak{C^{3}_{h}})=(3n+4)-(6n+5)+(3n+2)=1$; 
\item $r_{2}=\mathrm{rank}(\mathfrak{C^{2}_{h}})-\mathrm{rank}(\mathfrak{C^{3}_{h}})=(6n+5)-(3n+2)=3n+3$;  
\item $r_{3}=\mathrm{rank}(\mathfrak{C^{3}_{h}})=3n+2$. 
\end{itemize}
We need to show that $\mathrm{grade}(I_{r_{i}}(\mathfrak{C^{i}_{h}}))\geq i, 1 \leq i \leq 3$, where $I_{r_{i}}(\mathfrak{C^{i}_{h}})$ denotes 
the ideal generated by the $r_{i} \times r_{i}$ minors of the matrix 
$\mathfrak{C^{i}_{h}}$.
\medskip

\begin{itemize}
\item[(1)] We take following minors from $ \mathfrak{C_{h}^{2}}$:

\begin{eqnarray*}
\mathfrak{C^{[21]}_{h}} & := & [2 \,3\,\dots\, (3n+4)\,\vert\,1 \, \dots \,(3n+1) \,(6n+3)\, (6n+4)]\\
{} & = & -x_{1}^{r-n+5}(x_{2}x_{3}^{3}x_{1}^{r+n-1}-x_{1}^{r}x_{2}^{n+1}x_{3}x_{4}-x_{1}^{r+n}x_{4}^{3}+x_{2}^{n-1}x_{3}^{3n-6}x_{0}x_{4}^{r-3n+9})\\
\mathfrak{C^{[22]}_{h}} & := & [1 \,3\, 4 \,\dots \,3n+4\,\vert 2n\, (2n+2)\,\dots \,(3n-1) \,(3n+2) \,(4n+1) \,\dots \,(6n+1) \,(6n+4) \, (6n+5)]\\
{} & = & x_{3}^{3n-3}x_{4}(x_{2}^{n}x_{4}-x_{1}^{n-1}x_{3}^{2})(x_{3}^{3n-2}x_{4}^{r-3n+4}x_{0}-x_{1}^{r-n+2}x_{2}^{n+1})
\end{eqnarray*}
We see that any two irreducible components in the factorization of $\mathfrak{C^{[21]}_{h}},\mathfrak{C^{[22]}_{h}}$ are coprime. Therefore $\mathfrak{C^{[21]}_{h}},\mathfrak{C^{[22]}_{h}}$ forms a regular sequence of length $2$. Hence $\mathrm{grade}(I_{r_{2}}(\mathfrak{C^{2}_{h})})\geq 2$.
 \medskip
 
\item[(2)] We take the following minors from $\mathfrak{C^{3}_{h}}$:
\begin{eqnarray*}
\mathfrak{C^{[31]}_{h}} & := & [1 \,(2n+1) \,\dots\, (3n-2)\,  (4n+3)\, \dots \,5n \,(5n+1)\, (5n+3)\, (5n+4)\, \dots \\
{} & {} & (6n+1)(6n+4) \, \vert\, 1 \,2 \,\dots \,(3n+2)]\\
{} & = & x_{2}^{4n+1}-x_{1}^{2n-1}x_{2}^{2n}x_{3}x_{4};\\
\mathfrak{C^{[32]}_{h}} & := & [1\, 2 \,\dots\, 2n \,(3n+1)\, (3n+2)\, (3n+3)\, \dots \,4n \,\vert\, 1\, 2 \,\dots \,(3n+2)] =  x_{3}^{6n-2}-x_{2}x_{3}^{3n-1}x_{4}^{3n-2};\\
\mathfrak{C^{[33]}_{h}} & := & [3n\, (3n+1)\, \dots\,(6n+2)\, (6n+5)\,\vert\, 1\, 2 \,\dots\, (3n+2)] 
= x_{1}^{2r+4}.
\end{eqnarray*}

\noindent The leading terms of $\mathfrak{C^{[31]}_{h}},\mathfrak{C^{[32]}_{h}}$ 
and $\mathfrak{C^{[33]}_{h}}$, with respect to the negative degree reverse lexicographic 
ordering induced by $x_{1}<x_{2}<x_{3}<x_{4}<x_{0}$, are $x_{2}^{4n+1}, x_{3}^{6n-2}$ and 
$x_{1}^{2r+4}$, which are mutually coprime. Therefore, by Lemma 2.2 in \cite{Saha}, 
$\{\mathfrak{C^{[31]}_{h}},\mathfrak{C^{[32]}_{h}},\mathfrak{C^{[33]}_{h}}\}$ 
forms a regular sequence of length $3$ and hence $\mathrm{grade}(I_{r_{i}}(\mathfrak{C^{3}_{h}}))\geq 3$. 
This proves that the complex $\overline{\mathfrak{R}_{h}}$ is exact and its minimality 
follows from the fact that all entries of the 
matrices $\mathfrak{C_{h}^{i}}$,$1 \leq i \leq 3$, lie in the homogeneous maximal 
ideal $\langle x_{0},\dots,x_{4} \rangle$. Hence $\overline{\mathfrak{R}_{h}}$ 
is a minimal graded free resolution of $\overline{\mathfrak{P_{nr}}}$.
\end{itemize}

For the affine case we proceed as follows: We observe that $\mathfrak{\tilde{C}}_{h}$ is obtained by 
putting $x_{0}=1$ in the entries of $\mathfrak{C}_{h}$, hence it follows that 
$\mathfrak{\tilde{C}^{i}_{h}}\circ \mathfrak{\tilde{C}^{i+1}_{h}}=0$, $1 \leq i \leq 2$. 
After putting $x_{0}=1$ in the above minors, all the minors form a regular sequence and 
from Lemma \ref{Exact} we have the exactness of $\mathfrak{R}_{h}$. We also have that 
$\mathfrak{R}_{h}$ 
is a minimal graded free resolution of $\mathfrak{P_{nr}}$ since all the entries of the matrices 
$\mathfrak{\tilde{C}_{h}^{i}}$, $1 \leq i \leq 3$, lie in the homogeneous maximal 
ideal $\langle x_{1},\dots,x_{4} \rangle$. \qed

\section{Syzygies of the Projective closure of the Bresinsky curve}
In \cite{bre}, Bresinsky defined the family of curves $B_{h}$, for $ h\geq 2$, 
defined by the family of numerical semigroups $\Gamma ((2h-1)2h,(2h-1)(2h+1),2h(2h+1),2h(2h+1)+2h-1))$. 
Let $\frak{Q}_{h}$ denote the defining ideal of these curves. Bresinsky proved that 
the minimal number of generators $\mu(\frak{Q}_{h})$ is an unbounded function of 
$h$. Gr\"{o}bner basis and syzygies of the Bersinsky curve have been computed 
and it has been proved that all the Betti numbers are unbounded 
functions of $h$ in \cite{mssbetti}.  Here, we study the syzygies of 
the projecive closure of the Bresinsky curve and show that all the Betti 
numbers are not unbounded function of $h$. Note that, the last Betti number is $1$ 
for the projecive closure of the Bresinsky curve, but the projective curve is not   
Cohen-Macaulay.
 
\begin{theorem}\label{grobbresinsky}
Consider the degree reverse lexicographic monomial order induced by $ x_{1}>x_{2}>x_{3}> x_{4}$, in the polynomial ring $k[x_{1},x_{2},x_{3},x_{4}]$. For $h\geq 2$, let us consider the following polynomials:
\medskip

\noindent $p_{1}=x_{2}x_{3}-x_{1}x_{4}$\\[2mm]
$ p_{2}=x_{2}^{2h}-x_{3}^{2h-1}$\\[2mm]
$ p_{(3,j)}=x_{1}^{j+1}x_{3}^{2h-j}-x_{2}^{j}x_{4}^{2h-j},\, 0 \leq j \leq 2h-1 $\\[2mm]
$ p_{(4,i)}=x_{1}^{i+1}x_{2}^{2h-i}-x_{3}^{i-1}x_{4}^{2h-i},\, 1 \leq i \leq 2h $\\[2mm]
$ p_{5}=x_{3}^{4h}-x_{2}^{2h-1}x_{4}^{2h+1}$\\[2mm]
$ p_{(6,i)}=x_{1}^{2+i}x_{2}^{2h-2-i}x_{4}^{2+i}-x_{3}^{2h+1+i},\, 0 \leq i \leq 2h-3$\\[2mm]
$ p_{7}=x_{1}x_{2}^{2h-1}x_{4}-x_{3}^{2h} $\\[2mm]
$ p_{8}=x_{1}^{2h}x_{4}^{2h}-x_{3}^{4h-1}$.
\medskip

\noindent Let
$$\mathfrak{G}_{h}=\{p_{1},p_{2},p_{5},p_{7},p_{8}\}\cup\{p_{(3,j)}\mid 0 \leq j \leq 2h-1\}\cup \{p_{(4,i)}\mid 1 \leq i \leq 2h\}\cup \{p_{(6,i)}\mid 0 \leq i \leq 2h-3\}.$$ 
The set $\mathfrak{G}_{h}$ is a Gr\"{o}bner basis of $\frak{Q}_{h}$, with respect to the given monomial order. 
\end{theorem}

\proof We proceed by Buchberger's algorithm. We compute each $S$-polynomial and show that it reduces to zero upon division by $\mathfrak{G}_{h}$.
\medskip

\noindent $(1) \quad S(p_{1},p_{2})=x_{3}^{2h}-x_{1}x_{2}^{2h-1}x_{4}=-p_{7}$.
\medskip

\noindent $(2) \quad \mathrm{For}\,0 \leq j < 2h-1$,
\begin{align*}
S(p_{1},p_{(3,j)})=x_{2}^{j+1}x_{4}^{2h-l}-x_{1}^{2+j}x_{4}x_{3}^{2h-j-1}
&=-x_{4}(x_{1}^{j+1}x_{3}^{2h-j}-x_{2}^{j}x_{4}^{2h-j})\\
&=-x_{4}p_{(3,j+1)}.
\end{align*}

\begin{align*}
 S(p_{1},p_{(3,2h-1)})=x_{2}^{2h}x_{4}-x_{1}^{2h+1}x_{4}
&=x_{4}(x_{2}^{2h}-x_{3}^{2h-1})-x_{4}(x_{1}^{2h+1}-x_{3}^{2h-1})\\
&= x_{4}p_{2}-x_{4}p_{(4,2h)}.
\end{align*} 
\medskip

\noindent $(3)$ \quad For $ 1 \leq i < 2h$,
\begin{align*}
S(p_{1},p_{(4,i)})= x_{3}^{i}x_{4}^{2h-i}-x_{1}^{i+2}x_{2}^{2h-i-1}x_{4}
&=-x_{4}(x_{1}^{i+2}x_{2}^{2h-i-1}-x_{3}^{i}x_{4}^{2h-i-1})\\
&= -x_{4}p_{(4,i+1)},\\[2mm]
S (p_{1},p_{(4,2h)}) & = x_{2}x_{3}^{2h}-x_{1}^{2h+2}x_{4}.
\end{align*}
We have $\gcd(\mathrm{Lt}(p_{1}),\mathrm{Lt}(p_{(4,2h)}))=1$, therefore, 
$S (p_{1},p_{(4,2h)}) \longrightarrow_{\mathfrak{G}_{h}}0$. 
\medskip

\begin{align*}
(4) \quad S(p_{1},p_{5})=x_{2}^{2h}x_{4}^{2h+1}-x_{1}x_{3}^{4h-1}x_{4}&=(x_{2}^{2h}-x_{3}^{2h-1})x_{4}^{2h+1}-(x_{1}x_{3}^{2h}-x_{4}^{2h})x_{3}^{2h-1}x_{4}\\
&=p_{2}x_{4}^{2h+1}-p_{(3,0)}x_{3}^{2h-1}x_{4}.
\end{align*}

\begin{align*}
(5) \quad S (p_{1},p_{(6,i)})=x_{3}^{2h+2+i}-x_{1}^{3+i}x_{2}^{2h-3-i}x_{4}^{3+i}&=-p_{(6,i+1)}, 0 \leq i \leq 2h-4\\
&=-p_{8},  i=2h-3.
\end{align*}
\medskip

\noindent $(6) \quad S(p_{1},p_{7})=x_{3}^{2h+1}-x_{1}^{2}x_{2}^{2h-2}x_{4}^{2}=-p_{(6,0)}$.
\medskip

\begin{align*}
(7) \quad S(p_{1},p_{8})&=x_{2}x_{3}^{4h}-x_{1}^{2h+1}x_{4}^{2h+1}\\
&=x_{3}^{4h-1}(x_{2}x_{3}-x_{1}x_{4})+(x_{1}x_{3}^{2h}-x_{4}^{2h})x_{3}^{2h-1}x_{4}-x_{4}^{2h+1}(x_{1}^{2h+1}-x_{3}^{2h-1})\\
&=x_{3}^{4h-1}p_{1}+p_{(3,0)}x_{3}^{2h-1}x_{4}-x_{4}^{2h+1}p_{(4,2h)}.
\end{align*} 
\medskip

$(8) \quad S (p_{2},p_{(3,l)}) \longrightarrow_{\mathfrak{G}_{h}}0$, since 
$\gcd(\mathrm{Lt}(p_{2}),\mathrm{Lt}(p_{(3,l)}))=1$.
\medskip

\begin{align*}
(9) \quad S (p_{2},p_{(4,i)})&=x_{2}^{i}x_{3}^{i-1}x_{4}^{2h-1}-x_{1}^{i+1}x_{3}^{2h-1}\\
&=(x_{1}^{2}x_{3}^{2h-1}-x_{2}x_{4}^{2h-1})(-x_{1}^{i-1})+(x_{2}x_{3}-x_{1}x_{4})(\sum_{l=0}^{i-2}x_{1}^{l}x_{2}^{i-1-l}x_{3}^{i-2-l}x_{4}^{2h-i+l})\\
&=p_{(3,1)}(-x_{1}^{i-1})+p_{1}(\sum_{l=0}^{i-2}x_{1}^{l}x_{2}^{i-1-l}x_{3}^{i-2-l}x_{4}^{2h-i+l}),\quad 2 \leq i \leq 2h\\
S (p_{2},p_{(4,1)})&=-(x_{1}^{2}x_{3}^{2h-1}-x_{2}x_{4}^{2h-1})=-p_{(3,1)}.
\end{align*}
\medskip

\noindent $(10) \quad S(p_{2},p_{5})\longrightarrow_{\mathfrak{G}_{h}}0$, 
since $\gcd(\mathrm{Lt}(p_{2}),\mathrm{Lt}(p_{(5)}))=1$.  
\medskip

\begin{align*}
(11) \quad S(p_{2},p_{(6,i)})&=x_{2}^{2+i}x_{3}^{2h+1+i}-x_{1}^{2+i}x_{3}^{2h-1}x_{4}^{2+i}\\
&=-(x_{2}x_{3}-x_{1}x_{4})(\sum_{l=0}^{i+1}x_{1}^{l}x_{2}^{1+i-l}x_{3}^{2h+i-l}x_{4}^{l})\\
& = -p_{1}(\sum_{l=0}^{i+1}x_{1}^{l}x_{2}^{1+i-l}x_{3}^{2h+i-l}x_{4}^{l}),\quad  0 \leq i \leq 2h-3.
\end{align*}
\medskip

\noindent $(12) \quad S (p_{2},p_{7})=x_{2}x_{3}^{2h}-x_{1}x_{3}^{2h-1}x_{4}=(x_{2}x_{3}-x_{1}x_{4})x_{3}^{2h-1}=p_{1}x_{3}^{2h-1}$.
\medskip

\begin{align*}
(13) \quad S(p_{2},p_{8})=x_{2}^{2h}x_{3}^{4h-1}-x_{1}^{2h}x_{3}^{2h-1}x_{4}^{2h}&=(x_{2}x_{3}-x_{1}x_{4})(\sum_{l=0}^{2h-1}x_{1}^{l}x_{2}^{2h-1-l}x_{3}^{4h-2-l}x_{4}^{l})\\&=p_{1}(\sum_{l=0}^{2h-1}x_{1}^{l}x_{2}^{2h-1-l}x_{3}^{4h-2-l}x_{4}^{l}).
\end{align*}

\begin{align*}
(14) \quad S(p_{(3,i)},p_{(3,j)})&=x_{2}^{j}x_{3}^{j-i}x_{4}^{2h-j}-x_{1}^{j-i}x_{2}^{i}x_{4}^{2h-i}\\
&=(x_{2}x_{3}-x_{1}x_{4})(\sum_{r=0}^{j-i-1}x_{1}^{r}x_{2}^{j-1-r}x_{3}^{j-i-2-r}x_{4}^{2h-j+r})\\
&= p_{1}(\sum_{r=0}^{j-i-1}x_{1}^{r}x_{2}^{j-1-r}x_{3}^{j-i-2-r}x_{4}^{2h-j+r}).
\end{align*}
\medskip

\noindent $(15)$ \quad Let us consider three cases separately for $S(p_{(3,j)},p_{(4,i)})$.

\begin{enumerate}[(a)]
\item For $j = i$,
\begin{align*}
 S(p_{(3,j)},p_{(4,i)})
&=-x_{2}^{2h}x_{4}^{2h-i}+x_{3}^{2h-1}x_{4}^{2h-i}\\
&=-(x_{2}^{2h}-x_{3}^{2h-1})x_{4}^{2h-i}
=p_{2}(-x_{4}^{2h-i}).
\end{align*}

\item For $j < i$,
 \begin{align*}
S(p_{(3,j)},p_{(4,i)})& =x_{3}^{2h-j+i-1}x_{4}^{2h-i}-x_{1}^{i-j}x_{2}^{2h-i+j}x_{4}^{2h-j}\\
&= -x_{4}^{2h-i}(x_{1}^{i-j}x_{2}^{2h-i+j}x_{4}^{i-j}-x_{3}^{2h+i-j-1})
\\&=-x_{4}^{2h-i}p_{(6,i-j-2)}.
\end{align*}
\medskip

\item For $j > i$,
 \begin{align*}
S(p_{(3,j)},p_{(4,i)})&=x_{1}^{j-i}x_{3}^{2h+i-j-1}x_{4}^{2h-i}-x_{2}^{2h-i+j}x_{4}^{2h-i}\\
&=(x_{2}^{2h}-x_{3}^{2h-1})(-x_{2}^{j-i}x_{4}^{2h-j})\\
&-(x_{2}x_{3}-x_{1}x_{4})(\sum_{k=0}^{j-i-1}x_{1}^{k}x_{2}^{j-i-1-k}x_{3}^{2h-2-k}x_{4}^{2h-j+k})\\
&=p_{2}(-x_{2}^{j-i}x_{4}^{2h-j})-p_{1}(\sum_{k=0}^{j-i-1}x_{1}^{k}x_{2}^{j-i-1-k}x_{3}^{2h-2-k}x_{4}^{2h-j+k}).
\end{align*}
\end{enumerate}

\begin{align*}
(16) \quad S(p_{(3,j)},p_{5})&=x_{1}^{1+j}x_{2}^{2h-1}x_{4}^{2h+1}-x_{2}^{j}x_{3}^{2h+j}x_{4}^{2h-j}\\
&=(x_{1}x_{2}^{2h-1}x_{4}-x_{3}^{2h})x_{1}^{j}x_{4}^{2h}-(x_{2}x_{3}-x_{1}x_{4})(\sum_{k=0}^{j-1}x_{1}^{k}x_{2}^{j-1-k}x_{3}^{2h+j-1-k}x_{4}^{2h-j+k})\\
&=p_{7}x_{1}^{j}x_{4}^{2h}-p_{1}(\sum_{k=0}^{j-1}x_{1}^{k}x_{2}^{j-1-k}x_{3}^{2h+j-1-k}x_{4}^{2h-j+k}).
\end{align*} 
\medskip

\noindent $(17)$ \quad Let us consider four cases separately for $S(p_{(3,j)},p_{(6,i)})$.
\begin{enumerate}[(a)]
\item For $j=i$,
\begin{align*}
 S(p_{(3,j)},p_{(6,i)})&=x_{3}^{4h+1}-x_{1}x_{2}^{2h-2}x_{4}^{2h+2}\\
&=(x_{3}^{4h}-x_{2}^{2h-1}x_{4}^{2h+1})x_{3}+(x_{2}x_{3}-x_{1}x_{4})x_{2}^{2h-2}x_{4}^{2h+1}\\
&=p_{5}x_{3}+p_{1}x_{2}^{2h-2}x_{4}^{2h+1}. \\
\end{align*}

\item For $j=i+1$,
\begin{align*}
 S(p_{(3,j)},p_{(6,i)})&=x_{1}^{j-i-1}x_{3}^{4h+1+i-j}-x_{2}^{2h-2-i+j}x_{4}^{2h-j+i-2}\\&= p_{5}.
 \end{align*}
 
\item For $j>i+1$,
\begin{align*}
  S(p_{(3,j)},p_{(6,i)})&=x_{1}^{j-i-1}x_{3}^{4h+1+i-j}-x_{2}^{2h-2-i+j}x_{4}^{2h-j+i-2}\\
  &= (x_{1}x_{3}^{2h}-x_{4}^{2h})(x_{1}^{j-i-2}x_{3}^{2h+1+i-j})-(x_{2}^{2h}-x_{3}^{2h-1}(x_{2}^{j-i-2}x_{4}^{2h-j+i+2})\\&
  -(x_{2}x_{3}-x_{1}x_{4})(\sum_{k=0}^{j-i-3}x_{1}^{k}x_{2}^{j-i-3-k}x_{3}^{2h-2-k}x_{4}^{2h-j+i+2+k})
 \\
  &= p_{(3,0)}(x_{1}^{j-i-2}x_{3}^{2h+1+i-j})-p_{2}(x_{2}^{j-i-2}x_{4}^{2h-j+i+2})\\
&-p_{1}(\sum_{k=0}^{j-i-3}x_{1}^{k}x_{2}^{j-i-3-k}x_{3}^{2h-2-k}x_{4}^{2h-j+i+2+k}).
\end{align*}

\item For $j<i$,
\begin{align*}
S(p_{(3,j)},p_{(6,i)})&=(x_{3}^{4h}-x_{2}^{2h-1}x_{4}^{2h+1})x_{3}^{1+i-j}\\
&\,\,+(x_{2}x_{3}-x_{1}x_{4})(\sum_{k=0}^{i-j}(x_{1}^{k}x_{2}^{2h-2-k}x_{3}^{i-j-k}x_{4}^{2h+1+k})\\
&=p_{5}x_{3}^{1+i-j}+p_{1}(\sum_{k=0}^{i-j}(x_{1}^{k}x_{2}^{2h-2-k}x_{3}^{i-j-k}x_{4}^{2h+1+k}).
\end{align*}
\end{enumerate}
\medskip

\noindent $(18)$ \quad Let us consider two cases separately for $S(p_{(3,j)},p_{7})$.

\begin{enumerate}[(a)]
\item For $ 0 < j \leq 2h-1$,
\begin{align*}
S(p_{(3,j)},p_{7})&=-x_{2}^{2h+j-1}x_{4}^{2h-j+1}+x_{1}^{j}x_{3}^{4h-j}\\
&=(x_{1}x_{3}^{2h}-x_{4}^{2h})x_{1}^{j-1}x_{3}^{2h-l}+(x_{2}^{2h}-x_{3}^{2h-1})(-x_{2}^{j-1}x_{4}^{2h-j+1})\\&-(x_{2}x_{3}-x_{1}x_{4})(\sum_{k=0}^{j-2}x_{1}^{k}x_{2}^{j-2-k}x_{3}^{2h-2-k}x_{4}^{2h-j+1+k})\\
&=p_{(3,0)}x_{1}^{j-1}x_{3}^{2h-l}+p_{2}(-x_{2}^{j-1}x_{4}^{2h-j+1})
\\&-p_{1}(\sum_{k=0}^{j-2}x_{1}^{k}x_{2}^{j-2-k}x_{3}^{2h-2-k}x_{4}^{2h-j+1+k}).
\end{align*}

\item For $j=0$, we have $S(p_{(3,0)},p_{7})=x_{3}^{4h}-x_{2}^{2h-1}x_{4}^{2h+1}=p_{5}$.
\end{enumerate}

\allowdisplaybreaks
\begin{align*}
(19) \quad S(p_{(3,j)},p_{8})&=x_{3}^{6h-j-1}-x_{1}^{2h-j-1}x_{2}^{j}x_{4}^{4h-j}\\
&=(x_{3}^{4h}-x_{2}^{2h-1}x_{4}^{2h+1})x_{3}^{2h-j-1}+(x_{2}x_{3}-x_{1}x_{4})(\sum_{k=0}^{2h-j-2}x_{1}^{k}x_{2}^{2h-2-k}x_{3}^{2h-j-2-k}x_{4}^{2h+1+k})\\
&=p_{5}x_{3}^{2h-j-1}+p_{1}(\sum_{k=0}^{2h-j-2}x_{1}^{k}x_{2}^{2h-2-k}x_{3}^{2h-j-2-k}x_{4}^{2h+1+k}).
\end{align*}

\begin{align*}
(20) \quad S(p_{(4,i)},p_{(4,j)})&=x_{2}^{j-i}x_{3}^{j-i}x_{4}^{2h-j}-x_{1}^{j-i}x_{3}^{i-1}x_{4}^{2h-i}\\
&=(x_{2}x_{3}-x_{1}x_{4})(\sum_{k=o}^{j-i-1}x_{1}^{k}x_{2}^{j-i-1-k}x_{3}^{j-2-k}x_{4}^{2h-j+k})\\
&=p_{1}(\sum_{k=o}^{j-i-1}x_{1}^{k}x_{2}^{j-i-1-k}x_{3}^{j-2-k}x_{4}^{2h-j+k}).
\end{align*}
\medskip

\noindent $(21) \quad S(p_{(4,i)},p_{5})\longrightarrow _{\mathfrak{G}_{h}}0$, since 
$\gcd(\mathrm{Lt}(p_{(4,i)}),\mathrm{Lt}(p_{5}))=1$.
\medskip

\noindent $(22)$ \quad  Let us consider four cases separately for $S(p_{(4,i)},p_{(6,j)})$.

\begin{enumerate}[(a)]
\item For $ i=j$,
\begin{align*} 
 S(p_{(4,i)},p_{(6,j)})&=x_{2}^{2}x_{3}^{2h+1+i}-x_{1}x_{4}^{2+i}x_{3}^{i-1}x_{4}^{2h-i}\\
&=(x_{2}x_{3}-x_{1}x_{4})(x_{2}x_{3}^{2h+i}+x_{1}x_{3}^{2h+i-1}x_{4})+(x_{1}x_{3}^{2h}-x_{4}^{2h})x_{1}x_{3}^{i-1}x_{4}^2\\
&=p_{1}(x_{2}x_{3}^{2h+i}+x_{1}x_{3}^{2h+i-1}x_{4})+p_{(3,0)}x_{1}x_{3}^{i-1}x_{4}^2.
\end{align*}

\item For $ i<j$,
\begin{align*}
 S(p_{(4,i)},p_{(6,j)})&= x_{2}^{j+2-i}x_{3}^{2h+1+j}-x_{1}^{1+j-i}x_{3}^{i-1}x_{4}^{2h-i+j+2}\\
&=(x_{2}x_{3}-x_{1}x_{4})(\sum_{l=0}^{j-i+1}x_{1}^{l}x_{2}^{j+1-i-l}x_{3}^{2h+j-l}x_{4}^{l})\\
&+(x_{1}x_{3}^{2h}-x_{4}^{2h})x_{1}^{1+j-i}x_{3}^{i-1}x_{4}^{j-i+2}\\
&=p_{1}(\sum_{l=0}^{j-i+1}x_{1}^{l}x_{2}^{j+1-i-l}x_{3}^{2h+j-l}x_{4}^{l})+p_{(3,0)}x_{1}^{1+j-i}x_{3}^{i-1}x_{4}^{j-i+2}. \\
\end{align*}

\item For $ i \geq j+2$,
\begin{align*}
 S(p_{(4,i)},p_{(6,j)})&=x_{1}^{i-j-1}x_{3}^{2h+1+j}-x_{2}^{i-j-2}x_{3}^{i-1}x_{4}^{2h-i+j+2}\\
&=(x_{1}x_{3}^{2h}-x_{4}^{2h})x_{1}^{i-j-2}x_{3}^{1+j}-(x_{2}x_{3}\\&-x_{1}x_{4})(\sum_{l=0}^{i-j-3}x_{1}^{l}x_{2}^{i-j-3-l}x_{3}^{i-2-l}x_{4}^{2h-i+j+2+l})\\
&=p_{(3,0)}x_{1}^{i-j-2}x_{3}^{1+j}-p_{1}(\sum_{l=0}^{i-j-3}x_{1}^{l}x_{2}^{i-j-3-l}x_{3}^{i-2-l}x_{4}^{2h-i+j+2+l}).
\end{align*}

\item For $ i =j+1$,
\begin{align*} S(p_{(4,i)},p_{(6,j)})&=x_{2}x_{3}^{2h+1+j}-x_{3}^{j}x_{4}^{2h+1}=(x_{2}x_{3}-x_{1}x_{4})x_{3}^{2h+j}-(x_{1}x_{3}^{2h}-x_{4}^{2h})x_{3}^{j}x_{4}\\
&=p_{1}x_{3}^{2h+j}-p_{(3,0)}x_{3}^{j}x_{4}.
\end{align*}
\end{enumerate}

\begin{align*}
(23) \quad S(p_{(4,1)},p_{7})&=x_{3}^{2h}x_{1}-x_{4}^{2h}=p_{(3,0)}.\\
S(p_{(4,i)},p_{7})&=(x_{1}x_{3}^{2h}-x_{4}^{2h})x_{1}^{i-1}+(x_{2}x_{3}-x_{1}x_{4})(-\sum_{l=0}^{i-2}x_{1}^{l}x_{2}^{i-2-l}x_{3}^{i-2-l}x_{4}^{2h-i+1+l})\\
&=p_{(3,0)}x_{1}^{i-1}+p_{1}(-\sum_{l=0}^{i-2}x_{1}^{l}x_{2}^{i-2-l}x_{3}^{i-2-l}x_{4}^{2h-i+1+l}), \mathrm{for}\, 1 < i \leq 2h.
\end{align*}

\begin{align*}
(24) \quad S(p_{(4,i)},p_{8})&=x_{2}^{2h-i}x_{3}^{4h-1}-x_{1}^{2h-i-1}x_{3}^{i-1}x_{4}^{4h-i}\\
&=(x_{2}x_{3}-x_{1}x_{4})(\sum_{l=0}^{2h-i-1}x_{1}^{l}x_{2}^{2h-i-1-l}x_{3}^{4h-2-l}x_{4}^{l})+(x_{1}x_{3}^{2h}-x_{4}^{2h})x_{1}^{2h-i-1}x_{3}^{i-1}x_{4}^{2h-i}\\
&=p_{1}(\sum_{l=0}^{2h-i-1}x_{1}^{l}x_{2}^{2h-i-1-l}x_{3}^{4h-2-l}x_{4}^{l})+p_{(3,0)}x_{1}^{2h-i-1}x_{3}^{i-1}x_{4}^{2h-i},\mathrm{for}\,1 \leq i< 2h.\\
S(p_{(4,2h)},p_{8})&=x_{1}x_{3}^{4h-1}-x_{3}^{2h-1}x_{4}^{2h}\\
&=(x_{1}x_{3}^{2h}-x_{4}^{2h})x_{3}^{2h-1}\\
&=p_{(3,0)}x_{3}^{2h-1}.
\end{align*}

\begin{align*}
(25) \quad S(p_{(6,i)},p_{(6,j)})&=x_{2}^{j-i}x_{3}^{2h+1+j}-x_{1}^{j-i}x_{3}^{2h+1+i}x_{4}^{j-i}\\
&=(x_{2}x_{3}-x_{1}x_{4})(\sum_{l=0}^{j-i-1} x_{1}^{l}x_{2}^{j-i-1-l}x_{3}^{2h+j-l}x_{4}^{l})\\
&=p_{1}(\sum_{l=0}^{j-i-1} x_{1}^{l}x_{2}^{j-i-1-l}x_{3}^{2h+j-l}x_{4}^{l}),\,\mathrm{for} \, i < j.
\end{align*}

\begin{align*}
(26) \quad S(p_{(6,i}),p_{7})=x_{1}^{1+i}x_{4}^{1+i}x_{3}^{2h}-x_{2}^{1+i}x_{3}^{2h+1+i}&=-(x_{2}x_{3}-x_{1}x_{4})(\sum_{l=0}^{i}x_{1}^{l}x_{2}^{i-l}x_{3}^{2h+i-l}x_{4}^{l})
\\
&=-p_{1}(\sum_{l=0}^{i}x_{1}^{l}x_{2}^{i-l}x_{3}^{2h+i-l}x_{4}^{l}).
\end{align*}

\begin{align*}
(27) \quad S(p_{(6,i)},p_{8})&=x_{2}^{2h-2-i}x_{3}^{4h-1}-x_{1}^{2h-2-i}x_{4}^{2h-2-i}x_{3}^{2h+1+i}\\
&=(x_{2}x_{3}-x_{1}x_{4})(\sum_{l=0}^{2h-3-i}x_{1}^{l}x_{2}^{2h-3-i-l}x_{3}^{4h-2-l}x_{4}^{l})\\
&=p_{1}(\sum_{l=0}^{2h-3-i}x_{1}^{l}x_{2}^{2h-3-i-l}x_{3}^{4h-2-l}x_{4}^{l}).
\end{align*}

\begin{align*}
(28) \quad S(p_{7},p_{8})=x_{2}^{2h-1}x_{3}^{4h-1}-x_{1}^{2h-1}x_{3}^{2h}x_{4}^{2h-1}&=(x_{2}x_{3}-x_{1}x_{4})(\sum_{l=0}^{2h-2}x_{1}^{l}x_{2}^{2h-2-l}x_{3}^{4h-2-l}x_{4}^{l})\\
&=p_{1}(\sum_{l=0}^{2h-2}x_{1}^{l}x_{2}^{2h-2-l}x_{3}^{4h-2-l}x_{4}^{l}).
\end{align*}
\medskip

All the $S$-polynomials reduce to zero upon division by $\mathfrak{G}_{h}$. Therefore, 
$\mathfrak{G}_{h}$ is a Gr\"{o}bner basis of $\mathfrak{Q}_{h}$ with respect to the degree reverse lexicographic monomial order induced by $ x_{1}>x_{2}>x_{3}> x_{4}$ on $k[x_{1},x_{2},x_{3},x_{4}]$. 
\qed
\medskip

Given $g\in k[x_{1},x_{2},x_{3},x_{4}]$, let $g^{H}\in k[x_{0},x_{1},x_{2},x_{3},x_{4}]$ denote 
the homogenization of $g$, with respect to the indeterminates $x_{0}$.

\begin{lemma}\label{genprojbre-1}
For $h\geq 2$, let $\overline{\mathfrak{G}_{h}}:=\{g^{H}\mid g\in \mathfrak{G}_{h}\}$. 
Then $\overline{\mathfrak{G}_{h}}$ is a Gr\"{o}bner basis of the homogenized ideal 
$\overline{\mathfrak{Q}_{h}}$, with respect to the degree reverse lexicographic 
monomial order induced by $ x_{1}>x_{2}>x_{3}> x_{4}>x_{0}$ on $k[x_{0},x_{1},x_{2},x_{3},x_{4}]$.  
\end{lemma}

\proof Follows from Theorems \ref{grobbresinsky} and Theorem \ref{gbhom}.\qed

\begin{lemma}\label{genprojbre}
The ideal $\overline{\mathfrak{Q}_{h}}$, for $h\geq 2$, is generated by the following polynomials:
\medskip

\noindent $p_{1}^{H}=x_{2}x_{3}-x_{1}x_{4}$;\\[2mm]
$ p_{2}^{H}=x_{2}^{2h}-x_{0}x_{3}^{2h-1}$;\\[2mm]
$ p_{(3,j)}^{H}=x_{1}^{j+1}x_{3}^{2h-l}-x_{0}x_{2}^{j}x_{4}^{2h-j}$;  $ 0 \leq j \leq 2h-1 $;\\[2mm]
$ p_{(4,i)}^{H}=x_{1}^{i+1}x_{2}^{2h-i}-x_{0}^{2}x_{3}^{i-1}x_{4}^{2h-i}$;  $ 1 \leq i \leq 2h $;\\[2mm]
$ p_{5}^{H}=x_{3}^{4h}-x_{2}^{2h-1}x_{4}^{2h+1}$.
\end{lemma}

\proof  From Lemma \ref{genprojbre-1}, $\overline{\mathfrak{G}_{h}}=\{p_{1},p_{2},p_{5},p_{7},p_{8}\}\cup\{p_{(3,j)}\mid 0 \leq j \leq 2h-1\}\cup \{p_{(4,i)}\mid 1 \leq i \leq 2h\}\cup \{p_{(6,i)}\mid 0 \leq i \leq 2h-3\}$, is a Gr\"{o}bner basis of $\overline{\mathfrak{Q}_{h}}$. For $0 \leq i \leq 2h-3$, \begin{align*} 
p_{(6,i)}^{H}&=x_{1}^{2+i}x_{2}^{2h-2-i}x_{4}^{2+i}-x_{3}^{2h+1+i}\\
&=(x_{2}^{2h}-x_{0}x_{3}^{2h-1})x_{3}^{2+i}-(x_{2}x_{3}-x_{1}x_{4})\sum_{l=0}^{1+i}x_{1}^{l}x_{2}^{2h-1-l}x_{3}^{1+l}x_{4}^{l}\\
&=p_{2}^{H}x_{3}^{2+i}-p_{1}^{H}\sum_{l=0}^{1+i}x_{1}^{l}x_{2}^{2h-1-l}x_{3}^{1+l}x_{4}^{l}.
\end{align*} 

\begin{align*}
p_{7}^{H}=x_{1}x_{2}^{2h-1}x_{4}-x_{3}^{2h}x_{0}&=x_{3}(x_{2}^{2h}-x_{3}^{2h-1}x_{0})-x_{2}^{2h-1}(x_{2}x_{3}-x_{1}x_{4})\\
&=p_{2}x_{3}-p_{1}x_{2}^{2h-1}.
\end{align*}

\begin{align*}
p_{8}^{H}&=x_{1}^{2h}x_{4}^{2h}-x_{3}^{4h-1}x_{0}\\
&=x_{3}^{2h}(x_{2}^{2h}-x_{0}x_{3}^{2h-1})-(x_{2}x_{3}-x_{1}x_{4})(\sum_{l=0}^{2h-1}x_{1}^{l}x_{2}^{2h-1-l}x_{3}^{2h-1-l}x_{4}^{l})\\
&=x_{3}^{2h}p_{2}^{H}-p_{1}^{H}(\sum_{l=0}^{2h-1}x_{1}^{l}x_{2}^{2h-1-l}x_{3}^{2h-1-l}x_{4}^{l}).
\end{align*}
The set 
$\{p_{1}^{H},p_{2}^{H},p_{(3,j)}^{H},p_{(4,i)}^{H},p_{5}^{H}\mid 0 \leq j \leq 2h-1 ,1 \leq i \leq 2h \}$, for $h\geq 2$, generates the elements $p_{(6,i)}^{H}, p_{7}^{H}$, and $p_{8}^{H}$. 
Therefore, it follows that $\{p_{1}^{H},p_{2}^{H},p_{(3,j)}^{H},p_{(4,i)}^{H},p_{5}^{H}\mid 0 \leq j \leq 2h-1 ,1 \leq i \leq 2h \}$ generates $\overline{\mathfrak{Q}_{h}}$, for $h\geq 2$.  \qed

\begin{notations}\label{matrices}{\rm 
For writing the syzygies, we define the following matrices for $h\geq 2$:
\begin{enumerate}
\item $\mathfrak{B}^{1}_{h}:=\begin{pmatrix}
 \

p_{1}^{H} & p_{2}^{H} & p_{(3,0)}^{H} &\cdots & p_{(3,2h-1)}^{H}& p_{(4,1)}^{H} &\cdots & p_{(4,2h)}^{H}& p_{5}^{H}
\end{pmatrix}_{1\times (4h+3)}$
\medskip

\item $\mathfrak{B}^{2}_{h}= \begin{bmatrix}
& &  R^{B}_{1 \times 8h-2} \\
  A^{'}_{(2h+1)\times (2h+1)} & 0_{(2h+1) \times (2h-1)} & B^{'}_{(2h+1) \times (6h-1)} & 0_{(2h+1)\times (2h-1)} & D^{'}_{(2h+2)\times 6}\\
A^{''}_{(2h+1) \times (2h+1)} & B^{''}_{(2h+1) \times (2h-1)} & 0_{(2h+1) \times ((6h-1)} & C^{''}_{(2h+1) \times (2h-1)} & D^{''}_{(2h+1)\times 6}
\end{bmatrix}$ 
\medskip

where,
\medskip

\begin{itemize}
\item[$\bullet$]$R^{B}_{1 \times 8h-2}=\begin{bmatrix}
r_{1j}
\end{bmatrix}$, with
$r_{1,3+l}=x_{0}x_{2}^{l-1}x_{4}^{2h-l}, 1 \leq l \leq 2h-1. $\\[2mm]
$ r_{1,2h+1+l}=x_{0}^{2}x_{3}^{l-1}x_{4}^{2h-1-l}, 1 \leq l \leq 2h-1  $.\\[2mm]
$ r_{(1,4h+1+l)}=x_{1}^{l}x_{3}^{2h-l}, 1 \leq l \leq 2h-1  $.\\[2mm]
$ r_{(1,6h-1+l)}=x_{1}^{l+1}x_{2}^{2h-1-l}, 1 \leq l \leq 2h-1  $.\\[2mm]
$r_{1,j}=0, j \in \{1, \dots,8h-2\}\setminus 
\{3+l,2h+1+l,4h+1+l,6h-1+l\}, \text{for}\, 1 \leq l \leq 2h-1$.
\medskip

\item[$\bullet$]$ A^{'}_{(2h+1)\times (2h+1)}=\begin{bmatrix}
x_{0}x_{4} & x_{1}^{2} & 0 & 0 & 0 &\cdots &0\\
0 & 0 & -x_{1} & 0 & 0 & \cdots  &0 \\
0 & x_{0} & x_{3} & -x_{1} & 0 & \cdots & 0\\
0 & 0 &0 & x_{3} & \ddots & \ddots & 0 \\
\vdots & \vdots & \vdots & \ddots & \ddots & \ddots & 0 \\
0 & 0 & 0 & \cdots & \ddots & \ddots & -x_{1}  \\
0 & 0 & 0 & \cdots & \cdots & 0 & x_{3} 
\end{bmatrix}$
\medskip

\item[$\bullet$]$B^{'}_{(2h+1) \times (6h-1)}=
\begin{bmatrix}
0 & 0 & 0 &\cdots & 0\\
-x_{2} & 0 & 0 &\cdots & 0\\
x_{4} & -x_{2} & 0 &\cdots & 0\\
0& \ddots & \ddots & \ddots & 0 \\
\vdots & \ddots & \ddots& \ddots & 0 \\
0 & 0 & 0& x_{4} & -x_{2}\\ 
0 & 0 & 0& 0 & x_{4}\\ 
\end{bmatrix}$
\medskip

\item[$\bullet$]$D^{'}_{(2h+2)\times 6}=\begin{bmatrix}
x_{2}^{2h}-x_{0}x_{3}^{2h-1} & x_{1}x_{2}^{2h-1} & x_{0}x_{2}^{2h-1} & x_{2}^{2h-1}x_{4}^{2h} & x_{3}^{4h-1} & x_{2}^{2h-1}x_{3}^{2h}\\
-x_{2}x_{3}+x_{1}x_{4} & -x_{1}x_{3} & -x_{0}x_{3} & -x_{3}x_{4} & -x_{4}^{2h+1}& -x_{3}^{2h+1}\\
0 & -x_{0} & 0 & x_{3}^{2h} & x_{3}^{2h-1}x_{4} & x_{2}^{2h-1}x_{4}\\
0 & 0 & 0 & 0 & 0 & 0\\
\vdots & \vdots & \vdots & \vdots & \vdots & \vdots \\
0 & 0 &-x_{1} &0 &0&0\\

\end{bmatrix}$
\medskip

\item[$\bullet$]$A^{''}_{(2h+1) \times (2h+1)}=\begin{bmatrix}
0 & -x_{2} & 0 & \cdots & 0\\
0 & 0 & 0 & \ddots & 0\\
0 & \ddots & \ddots & \ddots & 0\\
0 & 0 & 0 & \cdots & 0\\
-x_{3} & 0 &0 & \cdots &0\\
0 & 0 & 0 & \cdots & 0\\
0 & 0 & 0 & \cdots & 0\\
\end{bmatrix}
$
\medskip

\item[$\bullet$]$B^{''}_{(2h+1) \times (2h-1)}=
\begin{bmatrix}
-x_{1} & 0 & \cdots & \cdots &0  \\
x_{2} & -x_{1} & \ddots & \ddots &0  \\
0 & x_{2} & \ddots & \ddots &0 \\
\vdots & \ddots & \ddots & \ddots & 0 \\
0 & 0 & \cdots & x_{2} & -x_{1} \\
0 & 0 & \cdots & 0 & x_{2}\\
0 & 0 & \cdots & 0 & 0\\
\end{bmatrix}
$
\medskip

\item[$\bullet$]$C^{''}_{(2h+1) \times (2h-1)}=
\begin{bmatrix}
-x_{3} & 0 & \cdots & \cdots &0  \\
x_{4} & -x_{3} & \ddots & \ddots &0  \\
0 & x_{4} & \ddots & \ddots &0 \\
\vdots & \ddots & \ddots & \ddots & 0 \\
0 & 0 & \cdots & x_{4} & -x_{3} \\
0 & 0 & \cdots & 0 & x_{4}\\
0 & 0 & \cdots & 0 & 0\\
\end{bmatrix},$
\medskip

\item[$\bullet$]$D^{''}_{(2h+1)\times 6}=\begin{bmatrix}
0 & x_{4} & 0 & \cdots &0& 0\\
0 & 0 & 0 & \cdots &0& 0 \\
\vdots & \vdots & \vdots & \vdots &\vdots& \vdots \\
0 & 0 & 0 & \cdots &0& 0 \\
0 & 0 &x_{3} & 0 &0 &0\\
0 & 0 &0 & -x_{1} &-x_{2} & -x_{0}\\
\end{bmatrix}$ 
\end{itemize}
\medskip

\item $\mathfrak{B}^{3}_{h}=\begin{bmatrix}
P_{(2h+1)\times (2h-1)} & Q_{(2h+1) \times (2h+2)} & \vline\\
 P^{'}_{(2h-1) \times (2h-1)} & Q^{'}_{(2h-1)\times (2h+2)} & \vline\\
  P^{''}_{(2h-1)\times (2h-1)} & Q^{''}_{(2h-1) \times (2h+2)} & \vline & R'_{(8h+4)\times 1} & R''_{(8h+4)\times 1}\\
  P^{'''}_{(2h-1) \times (2h-1)} &  Q^{'''}_{(2h-1) \times (2h+2)} & \vline\\
  0_{6 \times (2h-1)} & Q^{(iv)}_{6 \times (2h+2)} & \vline
\end{bmatrix}$, \, such that,
\begin{itemize}
\medskip

\item[$\bullet$]$P_{(2h+1)\times (2h-1)}=\begin{bmatrix}
0 & 0 &0 & \cdots & 0 & x_{3}\\
0 & 0  &0 & \cdots & 0 & 0\\
-x_{2} & 0 &0 & \cdots & 0 & 0\\
 x_{4}& -x_{2}& 0 & \cdots & 0 &0 \\
 0 & x_{4}  & -x_{2}& \ddots & 0 &0 \\
 \vdots & \ddots & \ddots & \ddots & 0 & 0\\
 0 & 0 & 0 & x_{4} & -x_{2} & 0 \\
  0 & 0 & 0 & 0 & x_{4} &-x_{2}  \\
\end{bmatrix}$
\medskip

\item[$\bullet$]$ P^{'}_{(2h-1) \times (2h-1)}=\begin{bmatrix}
0 &\cdots & 0 & 0\\
\vdots & \ddots &\vdots &\vdots\\
0  & \cdots & 0 & 0 \\
0  & \cdots & 0 & x_{1} \\
0  & \cdots & 0 & 0 \\
\end{bmatrix}$
\medskip

\item[$\bullet$]$ P^{''}_{(2h-1)\times (2h-1)} =\begin{bmatrix}
x_{1} & 0 &0 & \cdots & 0 &0\\
 -x_{3}& x_{1}& 0 & \cdots & 0 &0  \\
 0 & -x_{3}  & x_{1}& \ddots & 0 &0 \\
 \vdots & \ddots & \ddots & \ddots & 0 &0 \\
 0 & 0 & 0 & -x_{3} & x_{1} &0 \\
  0 & 0 & 0 & 0 & -x_{3} &0  \\

\end{bmatrix}$
\medskip

\item[$\bullet$]$ P^{'''}_{(2h-1) \times (2h-1)}=\begin{bmatrix}
0 &\cdots & 0 & 0\\
\vdots & \ddots &\vdots &\vdots\\
0  & \cdots & 0 & 0 \\
0  & \cdots & 0 & -x_{3} \\
\end{bmatrix}$
\medskip

\item[$\bullet$]$Q_{(2h+1) \times (2h+2)}=\begin{bmatrix}
0 & 0 &\cdots & 0& 0 & 0 & -x_{1} & 0\\
0 & 0 &\cdots & 0 &-x_{3} & -x_{4} &0 &0\\
0 & 0 &\cdots & 0 &x_{0} & 0 &0 & -x_{3}^{2h-1}x_{4} \\
0 & 0 &\cdots & 0 &0 & 0 &0 & 0 \\
\vdots & \vdots & \cdots & \vdots & \vdots & \vdots & \vdots & \vdots \\
0 & 0 &\cdots & 0& 0 & 0 & 0 & 0\\
\end{bmatrix}
$
\medskip

\item[$\bullet$]$ Q^{'}_{(2h-1)\times (2h+2)}=
\begin{bmatrix}
-x_{3} & 0 & 0 &\cdots & 0 & 0 & 0 &0 &0\\
x_{4} & -x_{3} & 0 & \ddots &0 & -x_{3} & -x_{4} & 0 & 0 \\
0 & x_{4} & -x_{3} & \ddots & 0 & x_{0} & 0 & 0 & 0\\
\vdots & \ddots & \ddots &  \ddots & 0 & 0 & 0 & 0 & 0\\
0 & 0 & 0 & x_{4}& -x_{3} & 0 & 0 & 0 & 0\\
0 & 0 & 0 & 0& x_{4} & 0 & 0 & x_{3} & 0\\
\end{bmatrix} $
\medskip

\item[$\bullet$]$ Q^{''}_{(2h-1) \times (2h+2)}=\begin{bmatrix}
0 & 0 & 0 &\cdots & 0 & 0 & 0 & x_{0} &0 & x_{3}^{2h}\\
0 & 0 & 0 &\cdots & 0 & 0 & 0 & 0 &0 & 0\\
\vdots & \vdots & \vdots & \vdots & \vdots & \vdots  & \vdots  & \vdots  &\vdots & \vdots\\
0 & 0 & 0 &\cdots & 0 & 0 & 0 & 0 &0 & 0\\
\end{bmatrix}$
\medskip

\item[$\bullet$]$ Q^{'''}_{(2h-1) \times (2h+2)}=\begin{bmatrix}
x_{1} & 0 & 0 &\cdots & 0 & x_{2} & 0 &0 &0\\
-x_{2} & x_{1} & 0 & \ddots &0 & 0 & 0 & 0 & 0 \\
0 & -x_{2} & x_{1} & \ddots & 0 & 0 & 0 & 0 & 0\\
\vdots & \ddots & \ddots &  \ddots & 0 & 0 & 0 & 0 & 0\\
0 & 0 & 0 & -x_{2}& x_{1} & 0 & 0 & 0 & 0\\
0 & 0 & 0 & 0& -x_{2} & 0 & 0 & 0 & 0\\
\end{bmatrix}$
\medskip

\item[$\bullet$]$Q^{(iv)}_{6 \times (2h+2)}=\begin{bmatrix}
0 & 0 & \cdots & 0 & x_{1} & x_{0} & -x_{4}^{2h}\\
0 & 0 & \cdots & -x_{1} & -x_{2} & 0 & 0\\
0 & 0 & \cdots & 0 & 0 & -x_{2} & 0\\
0 & 0 & \cdots & 0 & 0 & 0 & x_{2}\\
0 & 0 & \cdots & 0 & 0 & 0 & -x_{1}\\
0 & 0 & \cdots & 0 & 0 & 0 & 0
\end{bmatrix}$
\medskip

\item[$\bullet$]$R'_{(8h+4)\times 1}=\begin{bmatrix}
R'_{(l,1)}
\end{bmatrix} $ \,, such that,\\[2mm]
$R^{'}_{(8h-1,1)}=x_{3}^{2h}$,\\[2mm]
$R^{'} _{(4h+1+i,1)}=-x_{2}^{2h-1-i}x_{4}^{i+1}$,\\[2mm]
$R^{'}_{(8h,1)}=x_{3}^{2h-1}x_{4}$,\\[2mm]
$R^{'}_{(6h-1+j,1)}=x_{3}^{(2h-1-j}x_{4}^{j+1}$,\\[2mm]
$R^{'}_{(1,1)}=x_{4}^{2h},R^{'}_{(8h+3,1)}=x_{0}$,\\[2mm]
$R^{'}_{(8h+4,1)}=-x_{2}$,\\[2mm] 
$ R^{'}_{(l,1)}=0$, for $l\in [8h+4]\setminus\{8h-1,4h+1+i,8h,6h-1+j,1,8h+3,8h+4\}$,\\[2mm]
$0 \leq i \leq 2h-2$, \, $1 \leq j \leq 2h-2$.
\medskip

\item[$\bullet$]$R''_{(8h+4)\times 1}=\begin{bmatrix}
R''_{(l,1)}
\end{bmatrix} $\, , such that, \\[2mm]
 $ R^{''}_{(8h,1)}=x_{3}^{2h}$, \\[2mm]
 $R^{''}_{(3,1)}=-x_{2}^{2h-1}x_{4}$,\\[2mm]
 $R^{''}_{(4h+1+i,1)}=-x_{2}^{2h-1-i}x_{3}x_{4}^{i}$,\\[2mm]
 $R^{''}_{(6h-1+j,1)}=x_{3}^{2h-j}x_{4}^{j}$, \\[2mm]
 $R^{''}_{(1,1)}=x_{3}x_{4}^{2h-1}$,\\[2mm]
 $R^{''}_{(8h+2,1)}=x_{0}$,\\[2mm]
 $R^{''}_{(8h+4,1)}=-x_{1}$,\\[2mm]
 $R^{''}_{(l,1)}=0$, \, for $l\in [8h+4]\setminus\{8h,3,4h+1+i,6h-1+j,1,8h+2,8h+4\}$ for $1\leq i \leq 2h-2,1 \leq j \leq 2h-2$.
\end{itemize}
 \medskip
 
 \item $\mathfrak{B}^{4}_{h}:=\begin{pmatrix}
 \Gamma_{ij}
 \end{pmatrix}_{(4h+3)\times 1}$\, , such that, \, $\Gamma_{i1}, 1\leq i\leq 4h+3 $, \, are defined as follows:
 \medskip

\noindent $ \Gamma_{(4h-1,1)}=x_{3}^{2h}$,\\[2mm]
$ \Gamma_{(m,1)}=-x_{2}^{2h-m}x_{4}^{m}$, \, for $ 1 \leq m \leq 2h-1$,\\[2mm]
$ \Gamma_{(4h-2,1)}=-x_{3}^{(2h-1)}x_{4}$,\\[2mm]
$ \Gamma_{(2h-1+n,1)}=x_{3}^{2h-1-n}x_{4}^{n+1}$, \, for $ 1 \leq n \leq 2h-2$,\\[2mm] 
$ \Gamma_{(4h,1)}=-x_{4}^{2h}$,\\[2mm]
$ \Gamma_{(4h+1,1)}=-x_{0}$,\\[2mm]
$\Gamma_{(4h+2,1)}=-x_{1}$,\\[2mm]
$\Gamma_{(4h+3,1)}=x_{2}$,\\[2mm] 
$\Gamma_{(i,4h+1)}=0$, \, for $i\in \{1,\ldots,8h+4\}$ and \\[2mm]
$i\notin \{4h-1,m,4h-2,2h-1+n,4h,4h+1,4h+2,4h+3\}$,\\[2mm]
$1 \leq m \leq 2h-1, 1 \leq n \leq 2h-2$.

\end{enumerate}
}
\end{notations}

\begin{theorem}\label{resprojbresinsky}
For $h\geq 2$, a minimal graded free resolution of the homogenized ideal 
$\overline{\mathfrak{Q}_{h}}$, which is 
the defining ideal of the projective closure of the Bresinsky curve, over 
the polynomial ring $k[x_{0},x_{1},x_{2},x_{3},x_{4}]$, is given by
$$\mathbf{\mathfrak{B}_{h}}:\,0\longrightarrow R\stackrel{\mathfrak{B}^{4}_{h}}\longrightarrow R^{4h+3}\stackrel{\mathfrak{B}^{3}_{h}}\longrightarrow R^{8h+4}\stackrel{\mathfrak{B}^{2}_{h}}\longrightarrow R^{4h+3}\stackrel{\mathfrak{B}^{1}_{h}}\longrightarrow R\longrightarrow R/\overline{\mathfrak{Q}_{h}}\longrightarrow 0,$$
where the matrices $\mathfrak{B}^{i}, 1\leq i\leq 4$ have been defined above 
in notation \ref{matrices}.
\end{theorem}

\proof One can check that
$$\mathbf{\mathfrak{B}_{h}}:\, 0 \longrightarrow R\stackrel{\mathfrak{B}^{4}_{h}}\longrightarrow R^{4h+3}\stackrel{\mathfrak{B}^{3}_{h}}\longrightarrow R^{8h+4}\stackrel{\mathfrak{B}^{2}_{h}}\longrightarrow R^{4h+3}\stackrel{\mathfrak{B}^{1}_{h}}\longrightarrow R\longrightarrow R/\overline{\mathfrak{Q}_{h}}\longrightarrow 0$$
is a chain complex by calculating $\mathfrak{B}^{i}_{h}\mathfrak{B}^{i+1}_{h}=0,1\leq i\leq 3 $. 
To prove the exactness, we use Lemma \ref{Exact}. Let $r_{i}$ be the $i^{th}$ expected rank of $\mathbf{\mathfrak{B}_{h}^{i}} $. Then $r_{1}=(4h+3)-(8h+4)+(4h+3)-1=1$, $r_{2}=8h+4-(4h+3)+1=4h+2$, $r_{3}=4h+3-1=4h+2$, $r_{4}=1$. We need to show that $\mathrm{grade}(I_{r_{i}}(\mathfrak{B}_{h}^{i}))\geq i, 1\leq i\leq 4$, where $I_{r_{i}}(\mathfrak{B}_{h}^{i})$ denotes the ideal generated by 
the $r_{i}\times r_{i}$ minors of the matrix $\mathfrak{B}_{h}^{i}$. 
\medskip

We take $p_{1}^{H}\in I_{r_{1}}(\mathfrak{B}^{1})$, and we have 
$I_{r_{1}}(\mathfrak{B}^{1})\geq 1$.
\medskip

We take 
\begin{itemize}
\item $\mathfrak{L}^{[21]}_{h}:=[1\,3\,4\,\ldots \,(4h+3)|(4h+1)\,\ldots\, (8h-2)\,8h\,(8h+2)\,(8h+3)\,(8h+4)]$\\[2mm]
$\quad\quad =x_{4}^{4h-1}(-x_{2}^{2h}+x_{3}^{2h-1}x_{0})(x_{3}^{4h}-x_{2}^{2h-1}x_{4}^{2h+1})$.
\smallskip

\item $\mathfrak{L}^{[22]}_{h}:=[2\,\ldots \,(4h+3)|1 \,\,3 \,\,\ldots 4h\,\,(8h+1)\,(8h+2)\,\,(8h+4)]$\\[2mm]
$\quad\quad =-x_{1}^{2h-2}(-x_{2}x_{3}+x_{1}x_{4})(-x_{1}x_{2}^{2h-1}x_{4}+x_{3}^{2h}x_{0})(x_{0}x_{3}^{2h})$.
\end{itemize}
\medskip

\noindent The gcd of any two prime factors of $\mathfrak{L}^{[21]}_{h}$ and 
$\mathfrak{L}^{[22]}_{h}$ is $1$, therefore, $\mathfrak{L}^{[21]}_{h},\mathfrak{L}^{[22]}_{h}$ 
forms a regular sequence. Hence, $\mathrm{grade}(I_{r_{2}}(\mathfrak{B}^{2}_{h}))\geq 2 $.
\medskip

We take 
\begin{eqnarray*}
\bullet \quad \mathfrak{L}^{[31]}_{h} & := & [1\,(4h)\,\ldots\,(6h-3)\,\ldots \,(8h-3)\,(8h-1)(8h)\,(8h+3)\,(8h+4)\\
{} & {} & \quad \vert 1\,2\,\ldots\, (4h+1)\,(4h+3)]= x_{1}^{4h+2}
\end{eqnarray*}
\medskip

$\bullet \quad \mathfrak{L}^{[32]}_{h}:=[1\,\ldots\, 4h\,(8h+1)\,(8h+2)| 1\,\ldots\, (4h+2)\,]=x_{4}^{6h}(-x_{2}^{2h}+x_{3}^{2h-1}x_{0})$

\begin{align*}
\bullet \quad \mathfrak{L}^{[33]}_{h}&:=[3\,\ldots \,(2h+1)\,(4h+1)\,(6h)\,\ldots \,8h\,(8h+3)| 1\,\ldots\,(4h+2)]\\
&=x_{1}x_{2}^{6h-3}x_{4}-x_{0}(\sum_{i=0}^{2h-1} x_{1}^{i}x_{2}^{4h-2-i}x_{3}^{2h-i}x_{4}^{i})-x_{0}x_{1}^{2h-1}x_{2}^{2h-1}x_{3}x_{4}^{2h-1}\\
& \quad + x_{0}x_{1}^{2h-1}x_{2}^{2h-2}x_{3}x_{4}^{2h}+x_{0}^{2}x_{1}^{2h-2}x_{3}^{2h+1}x_{4}^{2h-2}.
\end{align*}
Consider the ideal $\langle\mathfrak{L}^{[31]}_{h},\mathfrak{L}^{[32]}_{h}\rangle$; its  
primary decomposition is 
$\langle x_{1}^{4h+2},x_{4}^{6h}\rangle\cap \langle x_{1}^{4h+2},(-x_{2}^{2h}+x_{3}^{2h-1}x_{0})\rangle$. Hence, the associated primes are $\langle x_{1},x_{4}\rangle$, $\langle x_{1},(-x_{2}^{2h}+x_{3}^{2h-1}x_{0})\rangle$. We observe that $\mathfrak{L}^{[33]}_{h}\notin\langle x_{1},x_{4}\rangle$ 
and $\mathfrak{L}^{[33]}_{h}\notin\langle x_{1},(-x_{2}^{2h}+x_{3}^{2h-1}x_{0})\rangle$. 
If $\mathfrak{L}^{[33]}_{h}\in\langle x_{1},x_{4}\rangle$, 
then $x_{0}x_{2}^{4h-2}x_{3}^{2h}\in \langle x_{1},x_{4}\rangle$, 
which is a contradiction, and if 
$\mathfrak{L}^{[33]}_{h}\in\langle x_{1},(-x_{2}^{2h}+x_{3}^{2h-1}x_{0})\rangle$, 
then again $x_{0}x_{2}^{4h-2}x_{3}^{2h}\in\langle x_{1},(-x_{2}^{2h}+x_{3}^{2h-1}x_{0})\rangle$, 
again a contradiction. Therefore, 
$\{ \mathfrak{L}^{[31]}_{h},\mathfrak{L}^{[32]}_{h},\mathfrak{L}^{[33]}_{h}\}$ forms 
a regular sequence; hence  $\mathrm{grade}(I_{r_{3}}(\mathfrak{B}^{3}_{h}))\geq 3 $. 
\medskip

We take $ \Gamma_{(4h,1)}=-x_{4}^{2h}$, \, $\Gamma_{(4h+1,1)}=-x_{0}$, 
\, $\Gamma_{(4h+2,1)}=-x_{1}$, \, $\Gamma_{(4h+3,1)}=x_{2}$. 
Therefore, $\mathrm{grade}(I_{r_{4}}(\mathfrak{B}^{4}_{h}))\geq 4$.
\medskip

We have proved that the complex $\mathbf{\mathfrak{B}_{h}}$ is exact and it 
is clear that all the entries of the matrices $\mathfrak{B}^{i}_{h}, 1\leq i\leq 4$, 
belong to the homogeneous maximal ideal $\langle x_{0},\ldots,x_{4}\rangle$. 
Hence, $\mathbf{\mathfrak{B}_{h}}$ is a minimal free resolution of 
the defining ideal of the projective closure of the Bresinsky curve $\overline{\mathfrak{Q}_{h}}$. \qed

\section{Syzygies of the projective closure of the Arslan Curve}
In \cite{a} Arslan gave a necessary and sufficient criterion for the 
Cohen-macaulayness of the tangent cone of a monomial curve at origin, 
using Gr\"{o}bner basis. He introduced the following family of curves $\mathfrak{A}_{h}$, 
for $h\geq 2$, defined by the family of numerical semigroups $\Gamma (h(h+1),h(h+1)+1,(h+1)^2,(h+1)^2+1))$, 
which we call the \textit{Arslan curve}. In this section, 
we find the syzygies of the projective closure of the Arslan curve, i.e., 
$\overline{\mathfrak{A}_{h}}$ and show that all the Betti numbers are unbounded 
functions of $h$.

\begin{lemma}
Let us consider the following polynomials, for $h\geq 2$: 
\medskip

\noindent $ w=x_{2}x_{3}-x_{1}x_{4}$;\\[2mm]
$ g_{i}=x_{1}^{i}x_{3}^{h-i+1}-x_{2}^{i+1}x_{4}^{h-i},$ \, for \, $ 0 \leq i \leq h-1$;\\[2mm]
$ g_{h}=x_{2}^{h+1}-x_{1}^{h}x_{3}$;\\[2mm]
$ q_{j}=x_{1}^{j+1}x_{2}^{h-j}-x_{3}^{j}x_{4}^{h-j}x_{0},$ for $ 0 \leq j \leq h $.
\medskip

\noindent The set $\mathfrak{U}_{h}=\{ w,g_{i},q_{j}\mid 1\leq i,j\leq h\}$ 
is a Gr\"{o}bner basis of the defining ideal of the projective closure 
of the Arslan curve $\mathfrak{A}_{h}$, with respect to the degree reverse 
lexicographic monomial order induced by $x_{1}>x_{2}>x_{3}>x_{4}>x_{0}$ \, 
on the polynomial ring $k[x_{0},x_{1},x_{2},x_{3},x_{4}]$ .
\end{lemma}

\proof We use Proposition 5.2 in \cite{hs} and Theorem \ref{gbhom}.

\begin{notations}\label{arsmatrix}{\rm 
For $h\geq 2$, let $\mathfrak{I}_{h}:=\mathfrak{p}(h(h+1),h(h+1)+1,(h+1)^2,(h+1)^2+1)$ and  $\overline{\mathfrak{I}_{h}}:=\overline{\mathfrak{p}(h(h+1),h(h+1)+1,(h+1)^2,(h+1)^2+1)}$. 
\medskip

\noindent For writing the syzygies, we first define the following matrices for $h \geq 2$:
\begin{enumerate}
\item $\mathfrak{A}^{1}_{h}:=\begin{pmatrix}
w&g_{0}&...&g_{h}&q_{0}&...&q_{h}
\end{pmatrix}$
\medskip

\item $\mathfrak{A}^{2}_{h}:=\begin{bmatrix}
\hspace{5cm} R \\
 \hline \\
A^{2}_{(h+1\times (h+1)} & B^{2}_{h \times h} & C^{2}_{h \times h} & 0_{h\times h}\\
A^{2'}_{(h+1) \times h} & 0_{(h+1) \times h } & 0_{(h+1) \times h }& C^{2'}_{(h+1) \times h}
\end{bmatrix}$ 
where,
\begin{itemize}
\item[$\bullet$]$A^{2}_{(h+1)\times (h+1)}=
\begin{bmatrix}
x_{0} & 0  & \cdots &0\\
0 & 0 & \cdots & 0\\
\vdots & \vdots & \vdots & 0\\
x_{1} & 0 & \cdots & 0
\end{bmatrix}$
\medskip

\item[$\bullet$]$B^{2}_{(h+1) \times h}=
\begin{bmatrix}
-x_{1} & 0 & \cdots & \cdots & 0\\
x_{3} & -x_{1} & \ddots & \vdots &\vdots \\
0 & x_{3} & \ddots &0 &0\\
\vdots & \ddots & \ddots & -x_{1} & 0\\
0 & 0 & & x_{3} & -x_{1}\\
0 & 0 & \cdots & 0 & -x_{3}
\end{bmatrix}
$
\medskip

\item[$\bullet$]$C^{2}_{h \times h}=
\begin{bmatrix}
-x_{2} & 0 & \cdots & 0\\
x_{4} & -x_{2} & \ddots &0 \\
0 & x_{4} & \ddots & 0\\
\vdots & \ddots & \ddots & -x_{2}\\
0 & 0 &\cdots & x_{4}
\end{bmatrix}$
\medskip

\item[$\bullet$]$A^{2'}_{(h+1) \times h}=
\begin{bmatrix}
-x_{2} & -x_{1} &0 & 0 & \cdots & 0\\
0 &x_{2} & -x_{1} & 0 & \cdots & 0\\
0 &0 & x_{2} & \ddots & \ddots &\vdots \\
0 &0 & 0 & \ddots & \ddots & 0 \\
\vdots & \vdots & \vdots & \ddots & \ddots & -x_{1}\\
0 & 0 &\cdots & \cdots& 0 & x_{2}
\end{bmatrix}$
\medskip

\item[$\bullet$]$C^{2'}_{(h+1) \times h}=
\begin{bmatrix}
-x_{3} & 0 & \cdots &  0\\
x_{4} & -x_{3} & \ddots  &0 \\
0 & x_{4} & \ddots  &0\\
\vdots & \ddots & \ddots  & -x_{3}\\
0 & 0 & \cdots& x_{4} \\
\end{bmatrix}$
and 
$R=\begin{bmatrix}
r_{1,1}&r_{2,1}&\cdots & r_{1,4h+1}
\end{bmatrix}$,
\medskip

such that,\\
$r_{1,l}=0, l \in \{1,\dots,4h+1 \} \setminus \{2+i,h+2+j,2h+1,2h+2+j,3h+1,3h+2+i\}$,\\
$ 0 \leq i \leq h-1, 0 \leq j \leq h-2$;\\[2mm]
$r_{1,2+i}=x_{3}^{i}x_{4}^{h-i-1}, 0 \leq i \leq h-1$;\\[2mm]
$ r_{1,h+2+j}=x_{2}^{j+1}x_{4}^{h-j-1},0 \leq j \leq h-2$;\\[2mm]
$r_{1,2h+1}=x_{2}^{h}$;\\[2mm]
$r_{1,2h+2+j}=x_{1}^{j}x_{3}^{h-j},0 \leq j \leq h-2$;\\[2mm]
$r_{1,3h+1}=x_{1}^{h-1}x_{3}$;\\[2mm]
$r_{1,3h+2+i}=x_{1}^{i+1}x_{2}^{h-i-1}, 0 \leq i \leq h-1$.
\end{itemize}
\medskip

\allowdisplaybreaks
\item $\mathfrak{A}^{3}_{h}:=\begin{bmatrix}
 0_{(h+1)\times (h-1)} & B^{3}_{(h+1)\times (h-1)} & \vline\\
 A^{3}_{h \times (h-1)} & 0_{(h)\times (h-1)} &\vline\\
 A^{3'}_{h \times (h-1)} & 0_{(h)\times (h-1)} &\vline & C^{3}_{(4h+1)\times 1}\\
 0_{h \times (h-1)} & B^{3'}_{h \times (h-1)} & \vline \\
\end{bmatrix}$
 where, 
 
 $A^{3}_{h \times (h-1)}=
\begin{bmatrix}
-x_{2} & 0 & \cdots & 0\\
x_{4} & -x_{2} & \ddots &0 \\
0 & x_{4} & \ddots & 0\\
\vdots & \ddots & \ddots & -x_{2}\\
0 & 0 &\cdots & x_{4}
\end{bmatrix}
$
,
$A^{3'}_{h \times (h-1)}=
\begin{bmatrix}
x_{1} & 0 & \cdots & 0\\
-x_{3} & x_{1} & \ddots &0 \\
0 & -x_{3} & \ddots & 0\\
\vdots & \ddots & \ddots & x_{1}\\
0 & 0 & & -x_{3}
\end{bmatrix}
$

$ B^{3}_{(h+1)\times (h-1)}=
\begin{bmatrix}
0 & 0 & \cdots & 0\\
-x_{3} & 0 & \cdots & 0\\
x_{4} & -x_{3} & \ddots &0 \\
0 & x_{4} & \ddots & 0\\
\vdots & \ddots & \ddots & -x_{3}\\
0 & 0 & & x_{4}
\end{bmatrix} $
,
$B^{3'}_{h \times (h-1)}=
\begin{bmatrix}
x_{1} & 0 & \cdots & 0\\
-x_{2} & x_{1} & \ddots &0 \\
0 & -x_{2} & \ddots & 0\\
\vdots & \ddots & \ddots & x_{1}\\
0 & 0 &\cdots & -x_{2}
\end{bmatrix}$
\medskip

$C^{3}_{(4h+1)\times 1}=\begin{bmatrix}
C_{i,1}
\end{bmatrix}$,\\[2mm]
such that,\\[2mm]
$C_{(1,1)}=-x_{2}x_{3}+x_{1}x_{4}$, \, $C_{(2,1)}=-x_{2}x_{4}$, \, $C_{(h+1,1)}=x_{3}^{2}$, \, 
$C_{(h+2,1)}=x_{0}x_{4}$,\\[2mm] 
$C_{(2h+1,1)}=-x_{1}x_{2}$, \, $C_{(2h+2,1)}=-x_{3}x_{0}$,\\[2mm]
$C_{(3h+1,1)}=x_{1}^{2}$, \, $C_{(3h+2,1)}=x_{2}^{2}$, \, $C_{(4h+1,1)}=-x_{1}x_{3}$,\\[2mm]
$C_{(i,1)}=0$, \, $i \in \{1, \dots,(4h+1)\}\setminus \{1, 2, h+1, h+2, 2h+1, 2h+2, 3h+1, 3h+2, 4h+1\}$.

\end{enumerate}
}
\end{notations}

\begin{theorem}
For $h\geq 2$, a graded free minimal resolution of $\overline{\mathfrak{I}_{h}}$, 
the defining ideal of the projective closure of the Arslan curve $\mathfrak{A}_{h}$, is
$$\mathbf{\mathfrak{M}_{h}}:\,0\longrightarrow R^{2h-1}\stackrel{\mathfrak{A}^{3}_{h}}\longrightarrow R^{4h+1}\stackrel{\mathfrak{A}^{2}_{h}}\longrightarrow R^{2h+3}\stackrel{\mathfrak{A}^{1}_{h}}\longrightarrow R\longrightarrow R/\overline{\mathfrak{I}_{h}}\longrightarrow 0,$$
where the matrices $\mathfrak{A}^{i}_{h} $, $1\leq i\leq 3$, are defined in notation \ref{arsmatrix}. 

\end{theorem}
\proof We use Lemma \ref{Exact}. Let $r_{i}$ be the $i^{th}$ expected rank of $\mathbf{\mathfrak{A}_{h}^{i}} $. Then $r_{1}=(2h+3)-(4h+1)+(2h-1)-1=1$, $r_{2}=4h+1-(2h-1)=2h+2$, $r_{3}=2h-1$. We need to show that $\mathrm{grade}(I_{r_{i}}(\mathfrak{A}_{h}^{i}))\geq i, 1\leq i\leq 3$. 
\medskip

We take the following minors from $ \mathfrak{A}^{2}_{h}$:
\begin{itemize}
\item $ \mathfrak{D}_{h}^{[21]} := [1\,3 \,\cdots \,(2h+3)|1 \,\cdots \,(2h+2)]=(x_{2}^{h+1}-x_{1}^{h}x_{3})(-x_{3}^{h+1}+x_{2}x_{4}^{h})x_{3}^{h-1}$;
\medskip

\item $ \mathfrak{D}_{h}^{[22]}:= [2 \,\cdots \,(2h+3)|1\,2\,(2h+2)\, \cdots \,(4h+1)]=(x_{2}x_{3}-x_{1}x_{4})(-x_{1}x_{2}^{h}+x_{4}^{h}x_{0})x_{4}^{h-1}$.

\end{itemize}
\medskip

\noindent It is clear from the factorization that $\mathfrak{D}_{h}^{[21]}, \mathfrak{D}_{h}^{[22]}$ 
form a regular sequence; hence $\mathrm{grade}(I_{r_{2}}(\mathfrak{A}_{h}^{2}))\geq 2$.
\medskip

We take following minors from $ \mathfrak{A}^{3}_{h}$:
\begin{itemize}
\item $\mathfrak{D}_{h}^{[31]} := [2\,\cdots\, h \,(2h+2)\,\cdots\, (3h+1)|1\,\cdots \,(2h-1)]$\\
$=x_{1}^{h+1}x_{3}^{h-1}-x_{3}^{2h-1}x_{0}=x_{3}^{h-1}(x_{1}^{h+1}-x_{3}^{h}x_{0})$;
\medskip

\item $ \mathfrak{D}_{h}^{[32]} := [2 \,\cdots \,2h|1\, \cdots\, (2h-1) ]=x_{2}^{h-1}x_{3}^{h+1}-x_{2}^{2h}x_{4}^{2h}=x_{2}^{h-1}(x_{3}^{h+1}-x_{2}^{h+1}x_{4}^{2h})$;
\medskip

\item $\mathfrak{D}_{h}^{[33]} :=  [(h+2)\,\cdots\, (2h+1)\,(3h+2) \,\cdots\, 4h|1 \,\cdots\, (2h-1)]$\\
$=x_{1}^{h}x_{2}^{h}-x_{1}^{h-1}x_{4}^{h}x_{0}=x_{1}^{h-1}(x_{1}x_{2}^{h}-x_{4}^{h}x_{0})$.
\end{itemize}
It can be seen that the primary decomposition of the ideal 
$\langle  \mathfrak{D}_{h}^{[31]},\mathfrak{D}_{h}^{[32]}\rangle$ is 
$$ \langle x_{3}^{h-1}, x_{2}^{h-1}\rangle\cap\langle x_{3}^{h-1},x_{3}^{h+1}-x_{2}^{h+1}x_{4}^{2h}\rangle\cap\langle x_{1}^{h+1}-x_{3}^{h}x_{0},x_{2}^{h-1}\rangle\cap\langle x_{1}^{h+1}-x_{3}^{h}x_{0},x_{3}^{h+1}-x_{2}^{h+1}x_{4}^{2h}\rangle.$$
Therefore, the associated primes of the ideal $\langle  \mathfrak{D}_{h}^{[31]},\mathfrak{D}_{h}^{[32]}\rangle$ are,
$$\langle x_{3}, x_{2}\rangle, \, 
\langle x_{3},x_{3}^{h+1}-x_{2}^{h+1}x_{4}^{2h}\rangle, \, 
\langle x_{1}^{h+1}-x_{3}^{h}x_{0},x_{2}\rangle, \, 
\langle x_{1}^{h+1}-x_{3}^{h}x_{0},x_{3}^{h+1}-x_{2}^{h+1}x_{4}^{2h}\rangle.$$
It is easy to verify that the polynomial $\mathfrak{D}_{h}^{[33]}$ does not belong 
to any associated primes of the ideal 
$\langle  \mathfrak{D}_{h}^{[31]},\mathfrak{D}_{h}^{[32]}\rangle$. Therefore, 
$\{\mathfrak{D}_{h}^{[31]},\mathfrak{D}_{h}^{[32]},\mathfrak{D}_{h}^{[33]} \}$ 
forms a regular sequence. Hence $\mathrm{grade}(I_{r_{3}}(\mathfrak{A}_{h}^{3}))\geq 3$. 
Minimality of the resolution follows from the fact that all the 
entries of the matrices $\mathfrak{A}^{i}_{h}, 1\leq i\leq 3$ belong to 
the homogeneous maximal ideal $\langle x_{0},\ldots,x_{4}\rangle$, and hence 
$\mathfrak{M}_{h}$ is a minimal free resolution of $\overline{\mathfrak{I}_{h}}$. \qed

\bibliographystyle{amsalpha}

\end{document}